\documentclass[11pt,dvipsnames,reqno]{amsart}
\usepackage[T1]{fontenc}
\usepackage{a4wide}
\usepackage{yhmath,mathrsfs,amsthm,amsmath,amssymb,amsfonts,enumerate,lipsum, appendix,mathtools, bm}
\usepackage{bbold,comment}
\usepackage{pdfpages}
\usepackage{orcidlink}
\usepackage[mathscr]{euscript}
\usepackage{graphicx}
\usepackage{pgf,tikz}
\usepackage{tkz-euclide}
\usepackage{graphicx}
\usepackage{subcaption}
\usetikzlibrary{shapes.geometric}
\usetikzlibrary{arrows,hobby}
\usepackage[shortlabels]{enumitem}
\usepackage[english]{babel}
\DeclareMathOperator\supp{supp}

\allowdisplaybreaks

\usepackage[sort,numbers]{natbib}

\usepackage{hyperref}
\usepackage{esint}
\hypersetup{
	colorlinks=true,
	linkcolor=b,
	filecolor=b,      
	urlcolor=b,
	citecolor=ao,
}
\definecolor{b}{rgb}{0.0, 0.0, 2.8}
\definecolor{ao}{rgb}{1.0, 0.13, 0.32}
\usepackage[margin=0.85in, heightrounded]{geometry}
\usepackage{bigints}
\usepackage{esint}
\newtheorem{theorem}{Theorem}[section]
\newtheorem{lemma}[theorem]{Lemma}
\newtheorem{proposition}[theorem]{Proposition}

\theoremstyle{definition}

\newtheorem{remark}[theorem]{Remark}
\numberwithin{equation}{section}
\usepackage[hyperpageref]{backref}

\newcommand*\rd{\mathbb{R}^N}

\newcommand*\N{\mathbb{N}}
\newcommand{\al} {\alpha}

\newcommand{\pa} {\partial}

\newcommand{\de} {\delta}

\newcommand{\om} {\omega}
\newcommand{\Om} {\Omega}
\newcommand{\la} {\lambda}

\newcommand{\si} {\sigma}

\newcommand{\Gr} {\nabla}
\newcommand{\no} {\nonumber}
\newcommand{\noi} {\noindent}
\newcommand{\eps} {\varepsilon}
\newcommand{\var} {\varepsilon}
\newcommand{\ra} {\rightarrow}

\newcommand{\bee} {\begin{equation}}
	\newcommand{\eee} {\end{equation}}
\newcommand{\bea} {\begin{eqnarray}}
	\newcommand{\eea} {\end{eqnarray}}
\newcommand{\Bea} {\begin{eqnarray*}}
	\newcommand{\Eea} {\end{eqnarray*}}

\def\d{\,{\rm d}}
\def\dx{\,{\rm d}x}
\def\dy{\,{\rm d}y}

\def\ds{\,{\rm d}s}

\def\RN{{\mathbb R}^N}
\def\({{\Big(}}
\def\){{\Big)}}

\def\dx{\,{\rm d}x}

\def\dxy{\,{\rm d}x\,{\rm d}y}

\DeclarePairedDelimiter\abs{\lvert}{\rvert}%
\DeclarePairedDelimiter\norm{\lVert}{\rVert}%

\vspace{-1.5\baselineskip}

\let\tmp\phi \let\phi\varphi \let\varphi\tmp

\newcommand{\frap}{(-\Delta_p)^s}
\newcommand{\delp}{-\Delta_p}
\newcommand{\fra}{(-\Delta)^s}
\newcommand{\del}{-\Delta}

\newcommand{\pst}{p^\ast}

\newcommand{\R}{\mathbb{R}}
\newcommand{\RR}{\mathbb{R}^N}

\newcommand{\ov}{\overline}

\renewcommand{\AA}{\mathcal{A}}

\newcommand{\DD}{\mathcal{D}}

\newcommand{\WW}{\mathcal{W}}
\newcommand\restr[2]{{
  \left.\kern-\nulldelimiterspace 
  #1 
  \littletaller 
  \right|_{#2} 
  }}
\usepackage{natbib}
\newcommand{\littletaller}{\mathchoice{\vphantom{\big|}}{}{}{}}

\title[Mixed local nonlocal equations with critical nonlinearity on $\mathbb{R}^N$]{Mixed local nonlocal equations with critical nonlinearity on $\mathbb{R}^N$: Non-existence, Existence, and Multiplicity of positive solutions}
\author[N. Biswas, S. Chakraborty and P. Das]{Nirjan Biswas\,\orcidlink{0000-0002-3528-8388}$^{1}$, Souptik Chakraborty\,\orcidlink{0009-0004-1867-0560}$^{2, \dagger}$, \and Paramananda Das\,\orcidlink{0009-0009-2821-6675}$^{3}$}
\email[N. Biswas]{nirjan.biswas@acads.iiserpune.ac.in, nirjaniitm@gmail.com} 
\email[S. Chakraborty]{souptik25@tifrbng.res.in, soupchak9492@gmail.com}
\email[P. Das]{paramananda.das@students.iiserpune.ac.in, pd348225@gmail.com}
\thanks{$^\dagger$Corresponding author}
\begin{document}
\maketitle
\centerline{$^{1,3}$Department of Mathematics, Indian Institute of Science Education and Research Pune,}
\centerline{Dr. Homi Bhabha Road, Pune 411008, India}
\centerline{$^2$Tata Institute of Fundamental Research, Centre For Applicable Mathematics,}
\centerline{Sharada Nagar, Bengaluru 560065, India}
\begin{abstract}
We consider the following quasilinear critical problem involving the mixed local nonlocal operator:
\begin{equation}\label{main_prob_abstract_1}\tag{$\mathcal{P}_p$}
        -\Delta_p u+(-\Delta_p)^s u=|u|^{p^*-2}u+f\text{ in }\mathbb{R}^N, 
\end{equation}
where $s \in (0,1), p \in (1, \infty), N>p$, $p^*=\frac{Np}{N-p}$, and $f$ is a nonnegative functional in the dual space of the ambient solution space. If $f \equiv0$, then we show that \eqref{main_prob_abstract_1} does not admit any nontrivial weak solution. This phenomenon stands in contrast to the purely local and purely nonlocal cases. On the other hand, if $f$ is a nontrivial nonnegative functional, we establish the existence of a positive weak solution to \eqref{main_prob_abstract_1} provided $\|f\|$ is small. For this purpose we prove the concentration compactness principle for the mixed operator $-\Delta_p +(-\Delta_p)^s$ in $\mathbb{R}^N$. We also discuss the multiplicity of positive weak solutions to \eqref{main_prob_abstract_1}. 
\end{abstract}
\maketitle

\section{Introduction}
For $N\ge 3$ and $s \in (0,1)$, we first consider the following class of critical problems driven by the mixed local nonlocal operator:
\begin{equation}\label{main_prob}\tag{$\mathcal{P}_2$}
        \del u+\fra u=|u|^{2^*-2}u+f \text{ in }\RR,\;  u \in \WW_2,
\end{equation}
where $2^*=\frac{2N}{N-2}$ is the critical Sobolev exponent, $f \in \WW_2^*$ is a nonnegative functional, where $\WW_2^*$ is the dual space of $\WW_2$, to be defined later.
The fractional Laplace operator is defined as
\begin{align*}
    \fra u(x) = \tilde{C}_{N,s} \lim_{\var \rightarrow 0^{+}} \int_{\RR \backslash B(x, \var)} \frac{u(x)-u(y)}{|x-y|^{N+2s}}\, \dy, \; \text{for}~x \in \RR,
\end{align*}
where $\tilde{C}_{N,s}>0$ is the normalization constant,  $B(x, \var)$ denotes the ball of radius $\var$ with centre at $x \in \RR$, and $\Delta u =\text{div}(\Gr u)$ is the classical Laplace operator.
The Gagliardo seminorm $[\cdot]_{s,2}$ is defined as
\begin{align*}
    [u]_{s,2}:= \left( \iint_{\R^{2N}}\frac{(u(x)-u(y))^2}{|x-y|^{N+2s}}\dxy \right)^{\frac{1}{2}}.
\end{align*}
From \cite[Theorem 1.1]{BDVV2025}, the following inequality holds
\begin{equation}\label{ineq-12-1}
    C_{N,s}\|u\|_{L^{2^*}(\RR)}^{2}\leq \left( \|\nabla u\|_{L^2(\RR)}^2 + [u]_{s,2}^2 \right), \; \text{ for all } \, u \in \mathcal{C}^{\infty}_c(\RR).
\end{equation}
In view of \eqref{ineq-12-1}, the natural solution space for \eqref{main_prob} is defined as
\begin{align*}
    \WW_2 := \overline{\mathcal{C}_c^{\infty}(\RR)}^{\|\nabla \cdot\|_{L^2(\RR)}+[\cdot]_{s,2}}.
\end{align*}
Using Brasco et al. ~\cite[Theorem~3.1, Appendix~A, B]{Brasco2021}, we can show that $\WW_2$ has the following characterization:
\begin{align*}
    \WW_2 = \left\{ u \in L^{2^*}(\RR) : \left( \int_{\RR} \abs{\Gr u}^2 \dx \right)^{\frac{1}{2}} + [u]_{s,2} < \infty \right\},
\end{align*}
and it is endowed with the norm $\rho_2(u):=\left(\|\nabla u\|_{L^2(\RR)}^2+[u]_{s,2}^2\right)^{\frac{1}{2}}$. In view of \cite[Lemma A.1 and Lemma B.1]{Brasco2021}, we infer that $\mathcal{C}_c^{\infty}(\RR)$ is dense in $\WW_2$. We say $u\in \WW_2$ is a weak solution of \eqref{main_prob} if
\begin{align}\label{test-function}
  \int_{\RR}\nabla u\cdot\nabla v\dx+\AA_2(u,v)=\int_{\RR}|u|^{2^*-2}uv\dx+\prescript{}{\WW_2^*}{\langle}f,v{\rangle}_{\WW_2}, \; \text{for all }\, v \in \WW_2,
\end{align}
where $$\AA_2(u,v) :=\iint_{\R^{2N}}\frac{(u(x)-u(y))(v(x)-v(y))}{|x-y|^{N+2s}}\dx\dy.$$

In the purely local case, the best Sobolev constant is only attained on $\R^N$ by the so-called Aubin-Talenti bubbles~\cite{Aubin76, Talenti1976}. These extremal functions (up to translation and dilation) provide all possible positive solutions to the associated Euler–Lagrange equation on $\R^N$ (see~\cite{GNN1979, CGS1989}). It also shows that the critical exponent problem in the local case does not admit a ground state solution on any bounded domain. The same behaviour holds in the purely nonlocal case (see \cite{Cotsiolis2004, Chen2006}). Consequently, research in this area naturally splits into two directions: one investigates how the topology of the underlying domain influences the existence of high-energy solutions, as demonstrated in the pioneering work of Bahri and Coron~\cite{Bahri-Coron1988}; the other examines perturbations of the critical equation, aiming to identify conditions under which positive solutions arise. A landmark contribution in this second direction was given by Br\'ezis and Nirenberg~\cite{BN1989}, who proved that although the unperturbed problem on a bounded domain has no solution, the perturbed problem with a nonzero functional always admits a minimizer for the associated constrained minimization problem on the unit $L^{2^*}$-sphere. Subsequent developments, notably by Tarantello~\cite{Tarantello92}, extended these results by establishing multiplicity of positive solutions under suitable conditions on the nonzero functional.

Recently, Biagi et al. \cite{BDVV2025} proved that the best Sobolev constant for the mixed local nonlocal case is not attained even on $\R^N$. Due to the lack of scale invariance of the mixed operator, it is not immediate whether the critical exponent problem in the mixed case on $\R^N$ has any nontrivial solution. Indeed, in Section~\ref{non_exst_section} we have shown that \eqref{main_prob} does not possess any nontrivial solution when $f \equiv 0$. This non-existence result is fundamentally distinct from both purely local and purely nonlocal cases, where we know that
the only positive solutions (up to scaling and translations) to $$\del u = \abs{u}^{2^*-2}u, \text{ in } \RR, \text{ and } \fra u = \abs{u}^{2^*_s-2}u \text{ in } \RR,$$ (where $2^*_s = \frac{2N}{N-2s}$) are $S_N(1+\abs{x}^2)^{-\frac{N-2}{2}}$ (local bubble) and $S_{N,s}(1+\abs{x}^2)^{-\frac{N-2s}{2}}$ (nonlocal bubble) respectively.
Thus, we ask the same question for the mixed operator:

\begin{center}
  \textbf{Question:} Does there exist a positive solution to the perturbed equation~\eqref{main_prob}?
\end{center}

The aim of this paper is to answer this question in the affirmative.
In the first part of this paper, depending on $\norm{f}_{\WW_2^*}$ and the dimension $N$, we examine the existence and multiplicity of positive weak solutions to \eqref{main_prob}.
In the purely local setup, several authors have studied the multiplicity of positive weak solutions to the following subcritical problem:
\begin{align}\label{p-local}
   -\Delta u + u= g(x,u(x)) + f \text{ in } \RR, \; u \in H^1(\RR),
\end{align}
where $g$ is a Carathéodory function, $g$ satisfies Ambrosetti-Rabinowitz type growth condition, and subcritical growth condition $\abs{g(x,t)} \le a\abs{t}^{r-1}+b\abs{t}$ where $r \in(2,2^*), a>0, b \in [0,1)$, $f \in (H^1(\RR))^*$ with $\norm{f}_{(H^1(\RR))^*}$ being sufficiently small. For example, Jeanjean in \cite{Jeanjean} provided the existence of at least two positive solutions to \eqref{p-local}.

In \cite{Jeanjean}, the author introduced a new variational approach which overcomes the degenerate structure of the set of possible critical points. This method uses the Palais-Smale (PS) decomposition of the energy functional corresponding to \eqref{p-local}. In \cite{AT2000}, Adachi-Tanaka provided the existence of at least four positive solutions to \eqref{p-local} taking $g(x,u(x))=a(x)|u|^{r-1}u$ with $r \in (1,2^*-1)$ and the coefficient function $a\geq 1$. They obtained three positive solutions below the first breaking-down level of Palais-Smale sequences using a local minimization argument and Lusternik-Schnirelmann category theory, and they found a fourth positive solution with energy between the first and second breaking-down levels of the Palais-Smale criterion. In the purely nonlocal setup of \eqref{p-local}, Bhakta-Chakraborty-Ganguly \cite{BhChGa} showed the existence of at least three positive solutions by using the approach of \cite{AT2000}. In the purely nonlocal setup, Bhakta-Pucci in \cite{Bhakta-Pucci}
provided the existence of at least two positive solutions to
\begin{align}\label{p-nonlocal}
    (-\Delta)^s u = \abs{u}^{2^*_s -2}u + f \text{ in } \RR, \; u \in \mathcal{D}^{s,2}(\RR),
\end{align}
when $f \in (\mathcal{D}^{s,2}(\RR))^*$ (see \eqref{D1pDsp} for the definition of $\DD^{s,2}(\RR)$). For that, they use a similar technique introduced in \cite{Jeanjean}. To establish the multiplicity of positive solutions, the authors analyze the characterization of the (PS) sequence of the energy functional corresponding to \eqref{p-nonlocal}.
A central element in this characterization is the existence of positive solutions to the limiting problem $(-\Delta)^s u = \abs{u}^{2^*_s -2}u$ in $\RR$.  The (PS) decomposition for \eqref{main_prob} is not available in the literature. For a bounded domain $\Omega$, the authors in \cite{ChGuMaSr} recently studied the (PS) decomposition of the Brezis-Nirenberg problem $\del u+\fra u=\la u +|u|^{2^*-2}u$ in $\Omega$, $u=0$ in $\RR \setminus \Omega$, $u\in X_0(\Omega)$, where $X_0(\Omega):= \{ u \in H^1(\RR): \restr{u}{\Omega} \in H^1_0(\Omega), \text{ and } u=0 \text{ a.e. in } \RR \setminus \Omega\}$. Due to the compactness of $X_0(\Omega) \hookrightarrow H^s(\Omega)$, the local Aubin-Talenti bubble (up to translation and scaling) appears in the (PS) decomposition. However, the embedding $\WW_2 \hookrightarrow \mathcal{D}^{s,2}(\RR)$ is not compact because of the non-compactness of $\WW_2 \hookrightarrow L^{2_s^*}(\RR)$.

Recently, in \cite{DXEZ2024}, Dipierro et al. studied the following subcritical mixed local-nonlocal problem
\begin{align}\label{D}
    -\Delta u + \fra u + u= u^{r-1} \text{ in } \RR, \; u>0 \text{ in } \RR, \; u \in H^1(\RR),
\end{align}
where $r \in (1,2^*)$. The existence of positive solutions for \eqref{D} is proved, relying on some new regularity results. In \cite{BaDeQu}, Barrios et al. showed the existence of positive solutions to the problem $\del u+\fra u=u^{p-1}\text{ in }\RR$, when $p<2^*$ but close to $2^*$.

We consider the following energy functional
$$J_f(u)=\frac12\rho_2(u)^2-\frac{1}{2^*}\|u_+\|_{2^*}^{2^*}-\prescript{}{\WW_2^*}{\langle}f,u{\rangle}_{\WW_2}, \;\text{for all }\, u \in \WW_2,$$
where $f\gneqq 0$, i.e., $\prescript{}{\WW_2^*}{\langle}f,\phi\rangle_{\WW_2 }\geq 0$ whenever $\phi \geq 0$ in $\WW_2$. Observe that every critical point of $J_f$ is a nonnegative weak solution to \eqref{main_prob}.

\begin{theorem}\label{Main_Theorem}
   Let $N \ge 3, s \in (0,1)$ and $f \in \WW_2^*$ with $f\gneqq 0$. Then there exist constants $r_0(N)>0$ and $d_1(N)>0$ such that if $\|f\|_{\WW_2^*}< d_1$, then \eqref{main_prob} admits a positive weak solution $u_{f,0}$ with $\rho_2(u_{f,0})<r_0$, and $$J_f(u_{f,0})=\inf_{\{\rho_2(u)< r_0\}}J_f(u)<0.$$ Further, if $N<6-4s$, then \eqref{main_prob} admits a second positive weak solution $u_{f,1}$ with $J_f(u_{f,1}) > J_f(u_{f,0})$.
\end{theorem}
For the existence of the first positive solution $u_{f,0}$, we show that $J_f$ is strictly convex in $B_{r_0}$. As we see in the proof of the above theorem, the precise choices of $r_0,\,d_1$ are the following:
$$r_0\leq \left(\frac{S_2^{\frac{2^*}{2}}}{2^*-1}\right)^{\frac1{2^*-2}} \; \text{ and } \; d_1(N)\leq \left(\frac12-\frac{1}{2^*(2^*-1)}\right)r_0.$$
The general strategy to prove the existence of the second positive solution is:
\begin{itemize}
    \item Get the first solution $u_{f,0}$ as a local minimizer which is a perturbation of $0$.
    \item Obtain an energy estimate of the form
    \begin{align}\label{e-1}
        J_f(u_{f,0}+R\Psi)<J_f(u_{f,0})+\frac1NS_2^{\frac N2},
    \end{align}
    for every $R\geq0$ and some positive $\Psi\in\WW_2$, where $S_2$ is the best constant of the embedding $\WW_2 \hookrightarrow L^{2^*}(\RR)$.
    \item Prove the (PS) condition below the level $J_f(u_{f,0})+\frac1NS_2^{\frac N2}$ and use the Mountain pass geometry.
\end{itemize}
In our situation, we get the following energy estimate (see Proposition \ref{energy_prop}):
\begin{align}\label{e-2}
J_f(u_{f,0}+RU_\eps)\leq J_f(u_{f,0})+\frac1NS_2^{\frac{N}{2}}+C_1\eps^{2-2s}-C_2\eps^{\frac{N-2}{2}}+o(\eps^{\frac{N-2}{2}}),
\end{align}
where $R>0$, $U_\eps$ is the product of the Aubin-Talenti bubble for the Laplace operator with a cutoff function. To get \eqref{e-1} from \eqref{e-2}, it is required that $N<6-4s$.
\begin{remark}
    Recall that we already have $N\geq3$. Hence the condition $N<6-4s$ restricts $s\in(0,\frac34).$ Indeed, the conditions $N<6-4s$ and $N\ge3$ imply that $N\in\{3,4,5\}$, which corresponds to $s<\frac34$, $\frac12$, and $\frac14$ respectively.
\end{remark}

Next, we consider the nonlinear analogue of \eqref{main_prob}, namely
\begin{equation}\label{main_prob_p}\tag{$\mathcal{P}_p$}
        \del_p u+\frap u=|u|^{p^*-2}u+f \text{ in }\RR, \; u\in\WW_p,
\end{equation}
where $p \in (1, \infty), p^*=\frac{Np}{N-p}$ is the critical Sobolev exponent, $f\gneqq 0$, and $f \in \WW_p^*$, where $\WW_p^*$ is the dual space of $\WW_p$, where
\begin{align*}
    \WW_p := \overline{\mathcal{C}_c^{\infty}(\RR)}^{\|\nabla \cdot\|_{L^p(\RR)}+[\cdot]_{s,p}},
\end{align*}
with
$$[u]_{s,p} := \left( \iint_{\R^{2N}}\frac{|u(x)-u(y)|^p}{|x-y|^{N+sp}}\dxy \right)^{\frac{1}{p}}.$$
The definition of $\WW_p$ is also natural due to the non-linear analogue of \eqref{ineq-12-1}, namely
\begin{align*}
    C_{N,s,p}\|u\|_{L^{p^*}(\RR)}^{p}\leq \left( \|\nabla u\|_{L^p(\RR)}^p + [u]_{s,p}^p \right), \;\text{ for all }\, u \in \mathcal{C}^{\infty}_c(\RR).
\end{align*}
Again using~\cite[Theorem~3.1, Appendix~A, B]{Brasco2021}, we see that
\begin{align*}
    \WW_p = \left\{ u \in L^{p^*}(\RR) : \left( \int_{\RR} \abs{\Gr u}^p \dx \right)^{\frac{1}{p}} + [u]_{s,p} < \infty \right\},
\end{align*}
and it is endowed with the norm $\rho_p(u):= \left(\|\nabla u\|_{L^p(\RR)}^p+[u]_{s,p}^p\right)^{\frac{1}{p}}$.
A function $u\in \WW_p$ is said to be a weak solution of \eqref{main_prob_p}, if

\begin{align}\label{test-function-p}
  \int_{\RR}\abs{\nabla u}^{p-2}\nabla u\cdot\nabla v\dx+\AA_p(u,v) 
  =\int_{\RR}|u|^{p^*-2}uv\dx+\prescript{}{\WW_p^*}{\langle}f,v{\rangle}_{\WW_p}, \; \text{ for all }\, v \in \WW_p,
\end{align}
where $$\AA_p(u,v) :=\iint_{\R^{2N}}\frac{\abs{u(x)-u(y)}^{p-2}(u(x)-u(y))(v(x)-v(y))}{|x-y|^{N+sp}}\dx\dy.$$
We consider the following energy functional associated with \eqref{main_prob_p}:
\begin{align*}
   J_f(u)=\frac1p\rho_p(u)^p-\frac1{p^*}\|u\|_{p^*}^{p^*}-\prescript{}{\WW_p^*}{\langle}f,u{\rangle}_{\WW_p}, \; \forall \, u \in \WW_p,
\end{align*}
where $f\gneqq 0$, i.e., $\prescript{}{\WW_p^*}{\langle}f,\phi\rangle_{\WW_p }\geq 0$ whenever $\phi \geq 0$ in $\WW_p$.

\begin{theorem}\label{Main Theorem_1}
     Let $p\in (1, \infty), s \in (0,1), N>p$, and $f \in \WW_p^*$ with $f\gneqq 0$. Then there exist constants $r_0(N,p)>0$ and $d_2(N,p)>0$ such that if $\|f\|_{\WW_p^*}< d_2$, then \eqref{main_prob_p} admits a positive weak solution $v_{f,0}$ with $\rho_p(v_{f,0})<r_0$, and $$J_f(v_{f,0})=\inf_{\{\rho_p(v)< r_0\}}J_f(v)<0.$$
\end{theorem}
We have noticed that the way of proving the first solution for \eqref{main_prob} is not adaptable for $p \neq 2$ (see Remark \ref{pnot2}). To find the existence of a positive solution, we take a minimizing sequence of $J_f$ on $B_{r_0}$ and show that the sequence is in fact a Palais-Smale (PS) sequence. In order to establish the convergence of the (PS) sequence, we use the concentration-compactness principle associated with $\del_p+\frap$ (proved in Proposition \ref{CC}).

\begin{remark}
    Regarding the multiplicity of positive solutions of \eqref{main_prob_p} for $p>2$, the main concern is the following energy estimate
    \begin{align}\label{e-1-p}
     J_f(v_{f,0}+R\Psi)<J_f(v_{f,0})+\frac1NS_p^{\frac Np},
    \end{align}
    where $S_p$ is the best constant of the embedding $\WW_p \hookrightarrow L^{p^*}(\RR)$.
    Proceeding similarly as in the proof of \cite[Proposition 4.5]{BhBiDa}, we see that, when $\nabla v_{f,0}\in L_{\text{loc}}^{\infty}(\RR)$,
    \begin{align*}
        &J_{f}(v_{f,0}+R W_\eps) \leq J_f(v_{f,0})+\frac1NS_p^{\frac Np}+C_1\eps^{\zeta_1(1-s)}-C_2\eps^{\frac{N-p}{p}}+o(\eps^{\frac{N-p}{p}}),\quad2<p<3;\\
        &J_{f}(v_{f,0}+R W_\eps) \leq J_f(v_{f,0})+\frac1NS_p^{\frac Np}+C_1\eps^{2-2s}+C_2\eps^{\frac{2(N-p)}{p(p-1)}}-C_3\eps^{\frac{N-p}{p}}+o(\eps^{\frac{N-p}{p}}),\quad p\geq3,
    \end{align*}
where $p-1<\zeta_1<\min\{2,\frac{N(p-1)}{N-1}\}$ is an arbitrary number, $R>0$, and $W_\eps$ is the product of the Aubin-Talenti bubble for the $p$-Laplace operator with a cutoff function. As a consequence:
\begin{enumerate}[(a)]
    \item For $2<p<3$, if $\eps^{\zeta_1(1-s)}<\eps^{(p-1)(1-s)}$, then \eqref{e-1-p} holds provided $N<p+p(p-1)(1-s)$, and proceeding similarly as in Section \ref{2_prob_section}, we obtain the multiplicity of positive solutions.
    \item For $p\geq3$, because of $\frac{2(N-p)}{p(p-1)} \le \frac{N-p}p$, the energy estimate \eqref{e-1-p} may fail. In this case, we may not conclude the multiplicity of positive solutions.
\end{enumerate}
The distinction between the cases $p<3$ and $p\ge3$ follows directly from the following inequality (see \cite[(v) and (vi) of Lemma 2.5]{BhBiDa}):
\begin{align*}
    (1+a^2+2a\cos\theta)^{\frac t2} \leq \begin{cases}
        1+a^t+ta\cos\theta+Ca^{\zeta_1}, & \text{ if } 2 < p<3;\\
        1+a^t+ta\cos\theta+C(a^2+a^{p-1}), & \text{ if } p \ge 3,
    \end{cases}
\end{align*}
where $a\ge 0$, $\zeta_1\in[t-1,2]$ and $\theta\in[0,2\pi]$. The above inequality plays a key role in proving the energy estimate \eqref{e-1-p}.
\end{remark}

Now, we state the non-existence result for the limiting problem.
\begin{theorem}\label{non-existence theorem}
    Let $p\in (1, \infty), s \in (0,1),\, N>p$, and $f \equiv 0$. Then \eqref{main_prob_p} does not admit any nontrivial weak solution in $\WW_p$.
\end{theorem}

To obtain the non-existence result, we use Pohozaev's identity for $\del_p+\frap$ on $\RN$, recently established by Anthal-Garain in \cite{AG}. Their result requires certain regularity and integrability for a weak solution to \eqref{main_prob_p}. To apply their result, we have relaxed those assumptions (see the discussion in Remarks \ref{M1} and \ref{G1}).

The remainder of the paper is organized as follows. In the next section, we list a few useful results. Section \ref{2_prob_section} is dedicated to proving the multiplicity of positive solutions to \eqref{main_prob}. In Section \ref{p_prob_section}, we study \eqref{main_prob_p} and prove the concentration-compactness principle. This section contains the proof of Theorem \ref{Main Theorem_1}. In Section \ref{non_exst_section}, we establish the non-existence result.
\vspace{0.1cm}

\noi \textbf{Notation:}
We use the following notation and convention throughout the paper.
\begin{enumerate}[(i)]
    \item We denote the norm $\norm{\cdot}_{L^p(\RR)}$ as $\norm{\cdot}_p$.
    \item The Beppo-Levi spaces $\mathcal{D}^{1,p}(\RR)$ and $\mathcal{D}^{s,p}(\RR)$ are defined as
    \begin{align}\label{D1pDsp}
    &\mathcal{D}^{1,p}(\RR) :=\left\{ u \in L^{p^*}(\RR) : \norm{\nabla u}_p < \infty \right\}, \no \\
    &\mathcal{D}^{s,p}(\RR) := \left\{ u \in L^{p^*_s}(\RR) : [u]_{s,p} < \infty \right\},
\end{align}
where $p^* = \frac{Np}{N-p}$ and $p^*_s = \frac{Np}{N-sp}$ are the critical Sobolev exponents.
    \item $S_p$ is the best constant in the inequality
\begin{equation}\label{ineq-12}
    C\|u\|_{p^*}^{p}\leq \left( \|\nabla u\|_p^p + [u]_{s,p}^p \right), \; \forall \, u \in \WW_p.
\end{equation}
Following \cite[Theorem 1.1]{BDVV2025}, we note that $S_p$ is also the best constant of the classical Sobolev inequality (\cite[Theorem 2.4]{T2007}):
\begin{equation}\label{ineq-1}
    C\|u\|_{p^*}^{p}\leq\|\nabla u\|_p^p, \; \forall \, u \in \mathcal{D}^{1,p}(\RR).
\end{equation}
\item $X^*$ denotes the dual space of $X$ and $\overset{\ast}{\rightharpoonup}$ denotes the weak* convergence on $X^*$.
\item $f\gneqq 0$ implies $f$ is a nonnegative nontrivial functional in $\WW_p^*$, i.e., $\prescript{}{\WW_p^*}{\langle}f,\phi\rangle_{\WW_p}\geq 0$ whenever $\phi \geq 0$ in $\WW_p$.
\item The fractional derivative of $u$ is defined as $$|D^s u|^p(x):=\int_{\RR}\frac{|u(x+h)-u(x)|^p}{|h|^{N+sp}}\d h.$$For any $u\in \WW_p$, $[u]_{s,p}^p=\int_{\RR}|D^su|^p(x)\dx$.
\item We set $A_p(u)(x,y):=|u(x)-u(y)|^{p-2}(u(x)-u(y))$. Thus $$\AA_p(u,v)=\iint_{\R^{2N}}\frac{A_p(u)(x,y)(v(x)-v(y))}{|x-y|^{N+sp}}\dx\dy.$$
\item $u_+=\max\{u,0\}$ and $u_-=\max\{-u,0\}$ denote the positive part and negative part of $u$ respectively.
\item $B_r(x)$ denotes the open ball of radius $r$ with center at $x$, and $\overline{B}_r(x)$ denotes the closure of $B_r(x)$.
\item $\text{card} (A)$ denotes the cardinality of a set $A$.
\item $C$ denotes a generic positive constant.
\end{enumerate}

\section{Preliminaries}\label{Preliminary}
The following proposition proves a compactness result of $\WW_p$ into the space of locally integrable functions.
\begin{proposition}\label{compact embeddings_1}
Let $s \in (0,1), p \in (1, \infty)$, and $N>p$. Then $\WW_p \hookrightarrow L_{\text{loc}}^q(\RN)$ compactly for every $q \in (0, p^*)$.
\end{proposition}

\begin{proof}
We adapt the argument given in \cite[Proposition 2.1]{Bi}. Clearly  $\WW_p \hookrightarrow L_{\text{loc}}^q(\RN)$ for $q \in (0, p^*)$. First, we show that $\WW_p$ is compactly embedded into $L_{\text{loc}}^p(\RN)$. Let $(u_n)$ be a bounded sequence in $\WW_p$, and $K \subset \RN$ be a compact set. Then the sequence $(u_n \big|_K)$ is bounded in $L^p(K)$. Using \cite[Lemma A.1]{BrLiPa} for every $n \in \N$, we have
\begin{align}\label{estimate1}
    \sup_{\abs{h}>0} \int_{\RN} \frac{|u_n(x+h)-u_n(x)|^p}{\abs{h}^{sp}} \, \dx \le C(N) [u_n]_{s,p}^p.
\end{align}
The boundedness of $(u_n)$ in $\WW_p$ and \eqref{estimate1} confirms that for all $n \in \N$, $\int_{\RN} |u_n(x+h)-u_n(x)|^p \, dx \ra 0$ as $\abs{h} \ra 0$. Now by applying the Riesz-Fr\'{e}chet-Kolmogorov compactness theorem on $(u_n \big|_K) \subset L^p(K)$, we conclude $(u_n \big|_K)$ is relatively compact. Therefore, $\WW_p$ is compactly embedded into $L_{\text{loc}}^p(\RN)$. For $q \in [1,p)$, using $L_{\text{loc}}^p(\RN) \hookrightarrow L_{\text{loc}}^q(\RN)$ we directly get the compact embeddings of $\WW_p$ into $L_{\text{loc}}^q(\RN)$. Further, for $q \in (0,1)$, using $L_{\text{loc}}^1(\RR) \hookrightarrow L_{\text{loc}}^q(\RR)$, we get the compact embeddings of $\WW_p$ into $L_{\text{loc}}^q(\RN)$. Next, in the case when $q \in (p,p^*)$, we express $q=pt+ (1-t)p^*$, where $t = \frac{p^*-q}{p^*-p} \in (0,1)$. Using $\WW_p \hookrightarrow L_{\text{loc}}^q(\RN)$, the sequence $(u_n \big|_K)$ is bounded in $L^q(K)$. Applying the H\"{o}lder's inequality with the conjugate pair $(\frac{1}{t}, \frac{1}{1-t})$ and \eqref{estimate1} we obtain the following estimate for every $\abs{h}>0$ and $n \in \N$:
\begin{align*}
    &\int_{\RR} \abs{u_n(x+h)-u_n(x)}^q \, \dx \\& \le \left( \int_{\RR} |u_n(x+h)-u_n(x)|^p \, \dx \right)^t \left( \int_{\RR} \abs{u_n(x+h)-u_n(x)}^{p^*} \, \dx\right)^{1-t} \\
    & \le  C(N) \left( \abs{h}^s  [u_n]_{s,p} \right)^{pt}  \norm{u_n}_{p^*}^{p^*(1-t)} \le  C(N)\abs{h}^{spt} \rho_p(u_n)^q,
\end{align*}
where the last inequality follows using \eqref{ineq-12}. Again the boundedness of $(u_n)$ in $\WW_p$, and the Riesz-Fr\'{e}chet-Kolmogorov compactness theorem confirm a convergent subsequence of $(u_n \big|_K)$ in $L^q(K)$. Therefore, $\WW_p$ is compactly embedded into $L_{\text{loc}}^q(\RN)$ for $q \in (p,p^*)$ as well. This completes the proof.
\end{proof}
The following two lemmas are useful in proving the concentration-compactness principle. See \cite{BoSaSi} for their proofs.

\begin{lemma}\label{Bonder_1}
   Let $s \in (0,1)$ and $p \in (1, \infty)$. Let $ \phi \in W^{1,\infty}(\RR)$ be such that $\supp(\phi)\subset B_1(0)$. Then
    \begin{equation}
        |D^s \phi|^p(x)\leq C(N,s,p,\|\phi\|_{W^{1,\infty}(\RR)})\min\left\{1,|x|^{-(N+sp)}\right\}.
    \end{equation}
\end{lemma}
\begin{lemma}\label{Bonder_2}
    Let $s \in (0,1), p \in (1, \infty)$, and $N>sp$. Let $w\in L^{\infty}(\RR)$ be such that
    $$0\leq w(x)\leq C|x|^{-\al},$$ where $\al>sp$ and $C>0$. Then $\WW_p \Subset L^p(w;\RR)$.
\end{lemma}

\begin{remark}
  By Lemma \ref{Bonder_1}, for any $\phi\in W^{1,\infty}(\RR)$ with compact support, $|D^s\phi|^p$ satisfies the hypotheses of Lemma \ref{Bonder_2} and thus $\WW_p \Subset L^p(|D^s\phi|^p;\RR)$.
\end{remark}

We also require the following Leibniz formula: for every $\theta > 0$, there exists $C_{\theta}>0$ such that
\begin{align}\label{Leibniz}
    \int_{\RN} \abs{D^s(uv)}^p \dx \le (1+ \theta) \int_{\RN} \abs{D^su}^p \abs{v}^p \dx 
    + C_{\theta} \int_{\R^N} \abs{D^s v}^p \abs{u}^p \dx, \; \forall \, u, v \in \WW_p.
\end{align}

Next, we state a strong maximum principle for \eqref{main_prob_p}. Its proof follows using the Logarithmic energy estimate for $u$ (see \cite[Lemma 3.4]{PK2022}) and a similar set of arguments given in \cite[Proposition 3.7]{NR2025}.
\begin{proposition}[Strong Maximum Principle]\label{SMP}
   Let $s \in (0,1)$ and $p \in (1, \infty)$. Assume that $u \in \WW_p$ is a weak supersolution of the following equation:
   \begin{align*}
       \delp u+ \frap u = 0  \text{ in } \mathbb{R}^{N},
\end{align*}
and $u \ge 0$ a.e. in $\R^N$. Then either $u \equiv 0$ or $u>0$ a.e. in $\R^N$.
\end{proposition}
Now we state a comparison principle for
\begin{equation}\label{pure_f_prob}
        \delp u+\frap u=f\text{ in }\RR, \; u\in \WW_p,\, f\in \WW_p^*.
    \end{equation}
\begin{proposition}[Comparison Principle]\label{comparison}
    Let $s\in(0,1)$ and $p \in (1, \infty)$. Let $u,v\in\WW_p$ be respectively a positive subsolution and a supersolution of \eqref{pure_f_prob}. Then $u\leq v$ a.e. in $\RR$.
\end{proposition}
\begin{proof}
From the hypothesis, for every $\phi \in \WW_p$ with $\phi \ge 0$,  we get
\begin{align*}
\left< \delp u - (\delp v), \phi \right> + \left< \frap u - (\frap v), \phi \right> \le 0,
\end{align*}
where $\left< \cdot, \cdot \right>$ denotes the duality action. Let $\Omega = \text{supp}((u-v)_+)$. Taking $\phi =(u-v)_+ \in \WW_p$ as a test function,
\begin{align}\label{SCP-1}
    I_1+I_2 \no &:=\int_{\RR} (\abs{\Gr u}^{p-2} \Gr u - \abs{\Gr v}^{p-2} \Gr v)\cdot\Gr((u-v)_{+}) \dx \no \\
    &\quad+ \iint_{\R^{2N}} \frac{\left( A_p(u)(x,y) - A_p(v)(x,y) \right)\left( (u-v)_+(x) - (u-v)_+(y)\right)}{\abs{x-y}^{N+sp}} \dxy \le 0.
\end{align}
We use the following inequality repeatedly: for $p \in (1, \infty)$ and $X, Y \in \R^N$,
\begin{equation}\label{com-inq}
\begin{aligned}
\left(\abs{X}^{p-2}X-\abs{Y}^{p-2}Y\right)(X-Y) \ge 0, \text{ and } \left(\abs{X}^{p-2}X-\abs{Y}^{p-2}Y\right)(X-Y) = 0 \text{ iff } X=Y.
\end{aligned}
\end{equation}
By \eqref{com-inq},
\begin{align*}
I_1=\int_{\Omega} (\abs{\Gr u}^{p-2} \Gr u - \abs{\Gr v}^{p-2} \Gr v)\cdot\Gr((u-v)) \dx \ge 0.
\end{align*}
Next, we split
\begin{align*}
    &\iint_{\R^{2N}} \frac{\left( A_p(u)(x,y) - A_p(v)(x,y) \right)\left( (u-v)_+(x) - (u-v)_+(y)\right)}{\abs{x-y}^{N+sp}} \dxy \\
    &= \left( \iint_{\Omega \times \Omega} + 2\iint_{\Omega \times (\RR \setminus \Omega)}\right) \frac{\left( A_p(u)(x,y) - A_p(v)(x,y) \right)\left( (u-v)_+(x) - (u-v)_+(y)\right)}{\abs{x-y}^{N+sp}} \dxy.
\end{align*}
Using \eqref{com-inq},
\begin{align*}
    \iint_{\Omega \times \Omega} \frac{\left( A_p(u)(x,y) - A_p(v)(x,y) \right)\left( (u-v)(x) - (u-v)(y)\right)}{\abs{x-y}^{N+sp}} \dxy \ge 0,
\end{align*}
For $x \in \Omega, y \in \RR \setminus \Omega$, note that $u(x)-v(x)\geq0\geq u(y)-v(y)$. Thus, we get $u(x)-u(y)\geq v(x)-v(y)$. Since the map $\R\ni \xi\mapsto |\xi|^{p-2}\xi$ is strictly increasing,
\begin{align*}
    \iint_{\Omega \times (\RR \setminus \Omega)} \frac{\left( A_p(u)(x,y) - A_p(v)(x,y) \right)(u-v)(x)}{\abs{x-y}^{N+sp}} \dxy \geq 0.
\end{align*}
Combining all the estimates along with \eqref{SCP-1}, we obtain $I_1,I_2=0$.
Now using the second part of \eqref{com-inq}, we observe that $$\nabla (u-v)^+=0 \text{ a.e. in }\RR.$$
Finally recalling that $(u-v)^+\in L^{p^*}(\RR)$, we conclude $u \le v$ a.e. in $\RR$.
\end{proof}
Next, we prove the boundedness of the weak solutions to \eqref{main_prob_p} when $f\equiv0$. This result is a key step toward proving Theorem \ref{non-existence theorem}. We require the following lemma.
\begin{lemma}
    Let $\xi,\eta\in\R$, $p>1$, and $\si\geq1$. Denote $\xi_\pm=\max\{\pm\xi,0\}$ and $\xi_M=\min\{\xi,M\}$. There exists $C(\si,p)>0$ such that the following inequalities hold
\begin{align}
  {\rm (i)} \; \,  &|\xi-\eta|^{p-2}(\xi-\eta) \left(\xi_+ ((\xi_+)_M)^{\sigma-1} - \eta_+ ((\eta_+)_M)^{\sigma-1}\right) \no \\
  &\ge C(\si,p)\left|\xi_+(\xi_+)_M^{\frac{\si-1}{p}} - \eta_+(\eta_+)_M^{\frac{\si-1}{p}} \right|^p.\label{(i2)ineq}\\
  {\rm (ii)} \; \, &|\xi-\eta|^{p-2}(\xi-\eta) \left(\eta_- ((\eta_-)_M)^{\sigma-1} - \xi_- ((\xi_-)_M)^{\sigma-1}\right) \no \\
  &\ge C(\si,p)\left|\xi_-(\xi_-)_M^{\frac{\si-1}{p}} - \eta_-(\eta_-)_M^{\frac{\si-1}{p}} \right|^p.\label{(i4)ineq}
\end{align}
\end{lemma}
\begin{proof}
We only prove \eqref{(i2)ineq}, while \eqref{(i4)ineq} follows by a similar argument. Clearly \eqref{(i2)ineq} holds when $$\{\xi\leq0,\eta\in\R\}\cup\{\xi,\eta\geq M\}$$ with $C(\si,p)=1$. When $0<\xi,\eta< M$, using \cite[Lemma A.1(i1)]{NR2025} we see that
\begin{align*}
        &|\xi-\eta|^{p-2}(\xi-\eta) \left(\xi_+ ((\xi_+)_M)^{\sigma-1} - \eta_+ ((\eta_+)_M)^{\sigma-1}\right)\\
        &=|\xi-\eta|^{p-2}(\xi-\eta) \left(\xi^{\sigma} - \eta^{\sigma}\right)\geq\left(\frac{p}{\si+p-1}\right)^p\left|\xi^{\frac{\si+p-1}{p}}-\eta^{\frac{\si+p-1}{p}}\right|^p\\
        &=\left(\frac{p}{\si+p-1}\right)^p\left|\xi_+(\xi_+)_M^{\frac{\si-1}{p}} - \eta_+(\eta_+)_M^{\frac{\si-1}{p}} \right|^p.
\end{align*}
Whereas when $0<\xi<M\leq \eta$, we let $x=\frac{\xi}{M}<1,y=\frac{\eta}{M}\ge1,$ to get
\begin{align*}
        &|\xi-\eta|^{p-2}(\xi-\eta) \left(\xi_+ ((\xi_+)_M)^{\sigma-1} - \eta_+ ((\eta_+)_M)^{\sigma-1}\right)\\
        &=|\eta-\xi|^{p-2}(\eta-\xi) \left(\eta M^{\si-1}-\xi^{\sigma}\right)=M^{\si+p-1}(y-x)^{p-1}(y-x^\si),\\
        &\left|\xi_+(\xi_+)_M^{\frac{\si-1}{p}} - \eta_+(\eta_+)_M^{\frac{\si-1}{p}} \right|^p=\left|\xi^{\frac{\si+p-1}p}-\eta M^{\frac{\si-1}{p}}\right|^p=M^{\si+p-1}\left(y-x^{\frac{\si+p-1}{p}}\right)^p.
\end{align*}
So we just need to show that $$C(\si,p)=\inf_{x\in(0,1), y\in[1,\infty)}f(x,y):=\inf_{x\in(0,1), y\in[1,\infty)}\frac{(y-x)^{p-1}(y-x^\si)}{\left(y-x^{\frac{\si+p-1}{p}}\right)^p}>0.$$
Clearly $f(x,y)$ is positive when $x\in(0,1), y\in[1,\infty)$. When $x\to0^+$ or $y\to\infty$, $f(x,y)\to1$. Now suppose $y=1+\eps$, $x=1-\de$ where $\eps\ge0,0<\de<\frac12$ and $\eps,\de\to0^+$. Then by Bernoulli's inequality and using the concavity of $1-(1-x)^\si$ on $[0,\frac12]$, we get
\begin{align*}
    \frac{(y-x)^{p-1}(y-x^\si)}{\left(y-x^{\frac{\si+p-1}{p}}\right)^p} &\ge\frac{(\eps+\de)^{p-1}(\eps+2(1-2^{-\si})\de)}{\left(\eps+\de\frac{\si+p-1}{p}\right)^p} \\
    &=\frac{\left(\frac{\eps}{\de}+1\right)^{p-1}\left(\frac{\eps}{\de}+2(1-2^{-\si})\right)}{\left(\frac{\eps}{\de}+\frac{\si+p-1}{p}\right)^p}\\
    &\ge\inf_{t\in[0,\infty)}\frac{\left(t+1\right)^{p-1}\left(t+2(1-2^{-\si})\right)}{\left(t+\frac{\si+p-1}{p}\right)^p}.
\end{align*}
It's easy to check that $$c(\si,p)=\inf_{t\in[0,\infty)}\frac{\left(t+1\right)^{p-1}\left(t+2(1-2^{-\si})\right)}{\left(t+\frac{\si+p-1}{p}\right)^p}>0.$$
Thus $f(1-\de,1+\eps)\geq c(\si,p)$ uniformly with respect to $\eps,\de$. This completes the proof of \eqref{(i2)ineq}.
\end{proof}

\begin{proposition}[$L^{\infty}$-bound]\label{regularity}
    Let $s \in (0,1),\, p \in (1, \infty),\, N>p$, and $f\equiv0$. Let $u\in \WW_p$ be a weak solution to \eqref{main_prob_p}. Then $u\in L^\infty(\RR)$.
\end{proposition}

\begin{proof}
\textbf{Uniform bound of the positive part:} For $M > 0$, $\si\geq1$, we define $$(u_+)_M= \min\{u_+, M\}, \text{ and } w_M=u_+ ((u_+)_M)^{\sigma-1}.$$ Clearly $0\leq(u_+)_M \in L^{\infty}(\RR)$ and using the Lipschitz continuity of the map $\R\ni t\mapsto t_+\min\{t_+,M\}^{\si-1}$, we see that $w_M \in \mathcal{W}_p$. For any $a \ge 1$, we take $\sigma = (a-1)p +1$. Then
\begin{align*}
\int_{\RR} |\Gr u|^{p-2} \Gr u \cdot \Gr w_M
&= \sigma \int_{\{u_+ \le M\}} |\Gr u_+|^p u_+^{\sigma-1} + \int_{\{u_+ > M\}} M^{\sigma-1} |\Gr u_+|^p\\
&= \sigma \int_{\{ u_+ \le M \}} |\Gr u_+|^p u_+^{(a-1)p} + \int_{\{ u_+ > M \}} M^{(a-1)p} |\Gr u_+|^p \\
&:= \sigma I_{\le M} + I_{>M}.
\end{align*}
On the other hand,
\begin{align*}
   \int_{\RR} |\Gr (u_+((u_+)_M)^{a-1})|^p
   &= \int_{\{ u_+ \le M \}}  |\Gr (u_+^{a})|^p + \int_{\{ u_+ > M \}} M^{(a-1)p} |\Gr u_+|^p \\
   &= a^p \int_{\{ u_+ \le M \}} |\Gr u_+|^p u_+^{(a-1)p} + \int_{\{ u_+ > M \}} M^{(a-1)p} |\Gr u_+|^p \\
   &= a^p I_{\le M} + I_{>M}.
\end{align*}
By Bernoulli's inequality,  $a^p = (1+(a-1))^p \ge 1 + p(a-1) \ge \sigma$ and thus
\begin{align}\label{claim-I}
    a^p I_{\le M} + I_{>M} \le \frac{a^p}{\sigma} \left( \sigma I_{\le M} + I_{>M} \right).
\end{align}
Moreover, using \eqref{(i2)ineq}, we see that
\begin{align}\label{nonloc_pos}
       & \iint_{\R^{2N}} \frac{A_p(u)(x,y) \big(w_M(x) - w_M(y)\big)}{|x-y|^{N+sp}}\dxy \nonumber\\
       &\ge C(\si,p)
       \iint_{\R^{2N}} \frac{\left| \big(u_+(x)(u_+)_M(x)\big)^{\frac{\sigma-1}{p}}-u_+(y)\big((u_+)_M(y)\big)^{\frac{\sigma-1}{p}} \right|^p}{\abs{x-y}^{N+sp}} \dxy \ge 0.
\end{align}
Therefore, taking $w_M$ as a test function, using \eqref{claim-I} and \eqref{nonloc_pos}, we obtain
\begin{align*}
\int_{\RR} |\Gr (u_+((u_+)_M)^{a-1})|^p & \le \frac{a^p}{\sigma} \int_{\RR} \abs{\Gr u}^{p-2} \Gr u \cdot \Gr w_M \\
 & \le \frac{a^p}{\sigma} \int_{\RR} \abs{u}^{p^*-2}u w_M = \frac{a^p}{(a-1)p+1}\int_{\RR}u_+^{p^*-p} (u_+((u_+)_M)^{a-1})^p.
\end{align*}
Set $v_M:= u_+((u_+)_M)^{a-1}$. Then the classical Sobolev inequality yields
\begin{align}\label{ite-0}
     S_p \left( \int_{\R^N} v_M^{p^*}\right)^{\frac{p}{p^*}} \le \frac{a^p}{(a-1)p+1}\int_{\RR} u_+^{p^*-p} v_M^{p},
\end{align}
for every $a \ge 1$. For $a_1=\frac{p^*}{p}$, we show that $u_+\in L^{a_1p^*}(\RR).$ Using H\"{o}lder's inequality with the conjugate pair $(\frac{p^*}{p}, \frac{p^*}{p^*-p})$, we write
\begin{align*}
   \int_{\RR} u_+^{p^*-p} v_M^{p} & \le \left( \int_{\{ u_+ \le R \}} + \int_{\{ u_+ > R\}}\right) u_+^{p^*-p} v_M^{p} \\
   & \le R^{p^*-p} \int_{\RR} v_M^{p} + \left( \int_{\{ u_+ > R\}} u_+^{p^*} \right)^{\frac{p^*-p}{p^*}} \left( \int_{\RR} v_M^{p^*} \right)^{\frac{p}{p^*}}.
\end{align*}
Now, choosing $R>0$ large enough so that
$$\left(\int_{\{u_+>R\}}u_+^{p^*}\right)^{\frac{p^*-p}{p^*}}\leq \frac{S_p((a_1-1)p+1)}{2a_1^p}.$$
Hence from \eqref{ite-0}, we get
\begin{align*}
    S_p \left( \int_{\R^N} v_M^{p^*} \right)^{\frac{p}{p^*}}
    \le \frac{a_1^p}{(a_1-1)p+1} \left(  R^{p^*-p} \int_{\RR} v_M^{p} + \frac{S_p((a_1-1)p+1)}{2a_1^p}  \left( \int_{\RR} v_M^{p^*} \right)^{\frac{p}{p^*}} \right),
\end{align*}
which implies
\begin{align*}
   S_p \left( \int_{\R^N} \left(u_+((u_+)_M)^{a_1-1}\right)^{p^*} \right)^{\frac{p}{p^*}} & \le \frac{2a_1^p}{(a_1-1)p+1}R^{p^*-p} \int_{\RR} v_M^{p}  \\ & \le \frac{2a_1^p}{(a_1-1)p+1}R^{p^*-p} \int_{\RR} u_+^{a_1p},
\end{align*}
Taking $M\to\infty$, and applying Fatou's lemma, we conclude that
\begin{align}\label{finite2}
    \left(\int_{\RR}|u_+|^{a_1p^*}\dx\right)^{\frac p{p^*}}\leq C(N,p)\int_{\RR}|u|^{p^*}\dx<\infty.
\end{align}
Now suppose $a> a_1$. Taking $M\to\infty$ in \eqref{ite-0}, we see that
\begin{align*}
    S_p\left(\int_{\RR}u_+^{ap^*}\right)^{\frac p{p^*}}\leq \frac{a^p}{(a-1)p+1}\int_{\RR}u_+^{p^*+(a-1)p}.
\end{align*}
Since $(a-1)p+1\geq a$, the above inequality implies that
\begin{align}\label{ite-1}
    \left(\int_{\RR}u_+^{ap^*}\right)^{\frac 1{p^*(a-1)}}\leq \left(\frac{a^{p-1}}{S_p}\right)^{\frac1{p(a-1)}}\left(\int_{\RR}u_+^{p^*+(a-1)p}\right)^{\frac1{p(a-1)}}.
\end{align}
We define $a_{m+1}$, $m\geq1$ by $$p^*+a_{m+1}p-p=p^*a_m,$$
which implies $$a_{m+1}-1=\left(\frac{p^*}p\right)^m(a_1-1).$$
Define $C_m=\frac{a_m^{p-1}}{S_p}$ and $$A_m=\left(\int_{\RR}u_+^{a_mp^*}\dx\right)^{\frac 1{p^*(a_m-1)}}.$$
Thus \eqref{ite-1} transforms to $$A_{m+1}\leq C_{m+1}^{\frac1{p(a_{m+1}-1)}}A_{m}\leq \Pi_{k=2}^{m+1}C_k^{\frac1{p(a_k-1)}}A_1.$$
This is exactly \cite[Eq. (4.12)]{BhBiDa}. Thus, $A_m\leq CA_1$ for every $m$, with $C$ independent of $m$. Therefore, taking $m \ra \infty$, from \eqref{finite2} we obtain $u_+\in L^\infty(\RR).$

\noi \noi  \textbf{Uniform bound of the negative part:} For $M > 0$, we similarly define $(u^-)_M= \min\{u^-, M\}$ and $z_M=-u_- ((u_-)_M)^{\sigma-1}$, for any $\sigma\geq1$. Clearly, $z_M \in \mathcal{W}_p$. Further, applying \eqref{(i4)ineq}, we also get
\begin{align*}
    \iint_{\R^{2N}} \frac{A_p(u)(x,y) \big(z_M(x) - z_M(y)\big)}{|x-y|^{N+sp}}\dxy \ge 0.
\end{align*}
Therefore, taking $z_M$ as a test function and using a similar set of iteration arguments, we get  $u_-\in L^\infty(\RR).$
Therefore, $u \in L^{\infty}(\RR)$.
\end{proof}

\section{Multiplicity of positive solutions for non-homogeneous problem}\label{2_prob_section}
Recall that
$$J_{f}(u)=\frac12\rho_2(u)^2-\frac{1}{2^*}\|u_+\|_{2^*}^{2^*}-\prescript{}{\WW_2^*}{\langle}f,u{\rangle}_{\WW_2}, \; \forall \, u \in \WW_2,$$ and observe $J_f \in \mathcal{C}^2(\WW_2;\R)$. Since $f\gneqq0$, by testing against $u_-$, we see that every critical point of $J_f$ is a nonnegative weak solution of \eqref{main_prob}. In the following proposition, we give the existence of a positive solution to \eqref{main_prob} that lies in a ball.
\begin{proposition}\label{first_sol}
    Let $s\in(0,1)$ and $f\gneqq0$. Then there exist $r_0>0$ and $d_1(r_0,N)>0$ such that $J_f$ is strictly convex in $B_{r_0}=\{u\in \WW_2:\rho_2(u)<r_0\}$. If $\|f\|_{\WW_2^*}< d_1$, then \begin{equation}\label{ex-1}
\delta_0=\inf_{\substack{\rho_2(u)=r_0,\\u\in \WW_2}}J_f(u)>0.
        \end{equation}
        Further, $J_f$ has a unique critical point $u_{f,0}\in B_{r_0}$ such that \begin{equation}\label{ex-2}
            0>c_{f,0}=J_f(u_{f,0})=\inf_{u\in B_{r_0}} J_f(u).
        \end{equation}
        Moreover, $u_{f,0}>0$ a.e. in $\RR$.
\end{proposition}
\begin{proof}
For $h \in \WW_2$,
    \begin{align*}
        J_f''(u)(h,h)&=\rho_2(h)^2-(2^*-1)\int_{\RR}u_+^{2^*-2}h^2\dx\\
        &\geq \rho_2(h)^2-(2^*-1)\|u_+\|_{2^*}^{2^*-2}\|h\|_{2^*}^2 \geq (1-(2^*-1)S_2^{-\frac{2^*}{2}}\rho_2(u)^{2^*-2})\rho_2(h)^2,
    \end{align*}
where the second inequality holds using the H\"{o}lder's inequality with the conjugate pair $(\frac{2^*}{2^*-2}, \frac{2^*}{2})$. We observe that  $$J_f''(u)(h,h) \ge 0, \text{ if } u \in B_{r} \text{ with } r\leq r_0\coloneqq \left(\frac{S_2^{\frac{2^*}{2}}}{2^*-1}\right)^{\frac1{2^*-2}}, $$ i.e., $J_f''(u)$ is positive definite and thus $J_f$ strictly convex in $B_{r_0}$. Next for $\rho_2(u)=r_0$, we see that
    \begin{align*}
        J_f(u)&=\frac12\rho_2(u)^2-\frac1{2^*}\|u_+\|_{2^*}^{2^*}-\prescript{}{\WW_2^*}{\langle}f,u{\rangle}_{\WW_2}\\
        &\geq \frac12 {r_0}^2-\frac1{2^*}S_2^{-\frac{2^*}2}{r_0}^{2^*}-\|f\|_{\WW_2^*}r_0\\
        &=\left(\frac12-\frac1{2^*}S_2^{-\frac{2^*}2}{r_0}^{2^*-2}\right)r_0^2-\|f\|_{\WW_2^*}r_0 =\left(\frac12-\frac1{2^*(2^*-1)}\right)r_0^2-\|f\|_{\WW_2^*}r_0.
    \end{align*}
    Hence there exists $d_1(r_0,N)>0$ such that if $\|f\|_{\WW_2^*}<d_1$, then \eqref{ex-1} holds. Now using the convexity of $J_f$ and \eqref{ex-1}, a minimizer $u_{f,0}\in B_{r_0}$ in \eqref{ex-2} exists and by the strict convexity, minimizer is unique. Since $J_f(0)=0$, $J_f(u_{f,0})=c_{f,0}\leq0$. Finally, since $f\neq0$, $u_{f,0}$ is nontrivial and by the strong maximum principle (Proposition \ref{SMP}), $u_{f,0}>0$ a.e. in $\RR$.
    For any $v\in B_{r_0}$, we observe that $$J_f(tv)\leq \frac12 t^2\rho_2(v)^2-t\prescript{}{\WW_2^*}{\langle}f,v{\rangle}_{\WW_2}.$$
    Since $f\gneqq0$, for $t>0$ small enough and for a nonnegative $v$, we get $J_f(tv)<0$ which infers $c_{f,0}=\inf_{B_{r_0}}J_f<0$.
    \end{proof}

\begin{remark}\label{pnot2}
Observe that the above argument is specific to the case $p=2$. Indeed, for $1<p<2$, the energy functional is not of class $\mathcal{C}^2$. On the other hand, for $p>2$, although the energy functional is of class $\mathcal{C}^2$, second derivative of the energy functional is considerably more involved and, in general, is not positive definite. Thus, the strict convexity of the energy functional on a small ball around the origin can not be proved analogously for $p\neq 2$.
\end{remark}

The following remark confirms that $u_{f,0}$ is locally bounded below by a positive quantity.

\begin{remark}\label{locally lower bounded}
    $u_{f,0}$ is a positive supersolution of the problem $$\del u+\fra u=f(x)\text{ in }\RR,$$ which has a unique positive solution $\ov{u}_f$, by Riesz representation theorem and the strong maximum principle (Proposition \ref{SMP}). Thus $\ov{u}_f$ is a positive supersolution of the problem $$\del u+\fra u=0\text{ in }\RR.$$ By the weak Harnack inequality \cite[Lemma 8.1]{PK2022} and comparison principle (Proposition \ref{comparison}), for any $r>0$, $x_0\in\RR$, and for some $Q\in(0,1)$,
    \begin{align}\label{1-1}
        u_{f,0}(x)\geq \ov{u}_f(x) \geq C\left(\fint_{B_{r}(x_0)}\ov{u}_{f}^Q\dx\right)^{\frac1Q}>0, \text{ a.e. }x\in B_r(x_0).
    \end{align}
    For any compact set $K\subset\RR$, we can cover it by finitely many balls $B_r(x_i)$, $1\leq i\leq n$. Thus, $$u_{f,0}\geq c_K:=\min_{1\leq i\leq n}\left(\fint_{B_{r}(x_i)}\ov{u}_{f}^Q\dx\right)^{\frac1Q}>0 \text{ a.e. in }K.$$
\end{remark}

The following energy estimates play a key role in finding the second positive solution.

\begin{proposition}[Energy Estimate]\label{energy_prop}
Let $s\in(0,1)$, $N<6-4s$ and $f\gneqq 0$. Let $u_{f,0}$ be as in the Proposition \ref{first_sol}. Then there exist $R_0>0$ and a positive function $\Psi\in \WW_2$ such that for every $R\geq R_0$,
\begin{equation}\label{energy estimate}
 \begin{split}
    &J_f(u_{f,0}+R\Psi)<J_f(u_{f,0}),\\
    &J_f(u_{f,0}+tR_0\Psi)<J_f(u_{f,0})+\frac1N S_2^{\frac{N}{2}},\quad\forall\,t\in[0,1].
\end{split}
\end{equation}
\end{proposition}
\begin{proof}
We follow the arguments of \cite{BhBiDa,BV2}. First we choose a Lebesgue point $y$ of $u_{f,0}$ in $\RR$, a cutoff function $\phi\in \mathcal{C}_c^{\infty}(B_{2r}(y))$ such that $0\leq\phi\leq1$ in $\RR$, $\phi\equiv1$ in $B_r(y)$ and $|\nabla\phi|\leq \frac{2}{r}.$ We then consider the family of scalings of the Talenti function:
    $$V_{\eps}(x)\coloneqq \frac{\left(N(N-2)\eps^2\right)^{\frac{N-2}{4}}}{\left(\eps^{2 }+|x-y|^{2}\right)^{\frac{N-2}{2}}},\quad \eps>0.$$
    We consider the following family of functions $U_\eps:=V_\eps\phi$. Using \cite[pg. 179]{St} and \cite[Lemma 5.3]{DFB}, as $\eps\to0$,
    \begin{equation}\label{e7.2}
        \begin{aligned}
       &\|\nabla U_\eps\|_2^2= S_2^{\frac N2}+O\left(\eps^{ (N-2)}\right), \\
        & \|U_\eps\|_{2^*}^{2^*}= S_2^{\frac N2}-O\left(\eps^{  N}\right),  \text{ and }\\
        &[U_\eps]_{s,2}^2=O\left(\eps^{ (N-2)}\right)+O\left(\eps^{ (2-2s)}\right).
    \end{aligned}
    \end{equation}
    Moreover, for $0< \eps< r$,
    \begin{align}\label{lb-2}
 \int_{B_r(y)}V_\eps^{2^*-1}\dx &=C\eps^{-\frac{N+2}{2}}\int_{B_r(y)}\frac{1}{\left(1+\left(\frac{|x-y|}{\eps}\right)^{2}\right)^{\frac{N+2}{2}}}\dx \no \\
 &=C\eps^{-\frac{N+2}{2}}\int_{0}^r\frac{\tau^{N-1}}{\left(1+(\frac{\tau}{\eps })^{2}\right)^{\frac{N+2}{2}}}\mathrm{d}\tau \no \\
 &=C\eps^{\frac{N-2}{2}}\int_0^{\frac{r}{\eps}}\frac{s^{N-1}}{(1+s^{2})^{\frac{N+2}{2}}}\ds \no \\
 &\geq C\eps^{ \frac{N-2}{2}}\int_0^{1}\frac{s^{N-1}}{(1+s^{2})^{\frac{N+2}{2}}}\ds=C\eps^{\frac{N-2}{2}}.
\end{align}
  For $0 \le t \le 1 \le R$, we set $w=u_{f,0}+tRU_{\eps}$. Then
    \begin{align}\label{(a.1)}
        J_f(w) & =\frac{\|\nabla u_{f,0}+tR\nabla U_{\eps}\|_2^2+[u_{f,0}+tRU_{\eps}]_{s,2}^2}{2} \no \\ & \quad -\frac{\|u_{f,0}+tRU_{\eps}\|_{2^*}^{2^*}}{2^*}-\prescript{}{\WW_2^*}{\langle}f,u_{f,0}+tRU_{\eps}{\rangle}_{\WW_2}.
    \end{align}
Expanding we obtain
\begin{align*}
J_f(w)&=\frac12\rho_2(u_{f,0})^2+\frac{(tR)^2}2\rho_2(U_\eps)^2+tR\left(\int_{\RR}\nabla u_{f,0}\cdot\nabla U_{\eps}\dx+\AA_2(u_{f,0},U_\eps)\right)\\&\quad-\frac{\|u_{f,0}+tRU_{\eps}\|_{2^*}^{2^*}}{2^*}-\prescript{}{\WW_2^*}{\langle}f,u_{f,0}+tRU_{\eps}{\rangle}_{\WW_2}.
\end{align*}
Since $u_{f,0}$ solves \eqref{main_prob}, we further have
\begin{align*}
J_f(w)&=\frac12\rho_2(u_{f,0})^2+\frac{(tR)^2}2\rho_2(U_\eps)^2 \\ & \quad -\frac1{2^*}\left(\|u_{f,0}+tRU_{\eps}\|_{2^*}^{2^*}-2^* tR\int_{\RR}u_{f,0}^{2^*-1}U_\eps\dx\right)-\prescript{}{\WW_2^*}{\langle}f,u_{f,0}{\rangle}_{\WW_2}\\
&=J_f(u_{f,0})+\frac{(tR)^2}2\rho_2(U_\eps)^2-\frac{(tR)^{2^*}}{2^*}\|U_\eps\|_{2^*}^{2^*} \\ & \quad -(tR)^{2^*-1}\int_{\RR}U_\eps^{2^*-1}u_{f,0}\dx-L_1,
\end{align*}
where
\begin{align*}
    L_1&:=\frac1{2^*}\int_{\RR}|u_{f,0}+tRU_{\eps}|^{2^*}-|u_{f,0}|^{2^*}-(tR)^{2^*}|U_{\eps}|^{2^*}\dx\\
    &\quad-\int_{\RR}tRu_{f,0}U_\eps(u_{f,0}^{2^*-2}+(tRU_\eps)^{2^*-2})\dx.
\end{align*}
Since $N<6-4s$, we have $2^*>3$. Now using the inequality $(1+a)^{2^*}\geq 1+a^{2^*}+2^*a+2^*a^{2^*-1}$ for every $a \ge 0$, we see that $L_1\geq0$. Using Remark \ref{locally lower bounded} and \eqref{lb-2}, we also have
\begin{align}\label{estimate-3}
    \int_{\RR}U_{\eps}^{2^*-1}u_{f,0}\dx\geq C\int_{B_r(y)}V_\eps^{2^*-1}\dx\geq C\eps^{ \frac{N-2}2}.
\end{align}
Combining all the above inequalities,
 \begin{equation}\label{(d1)}
     \begin{aligned}
         J_f(w)&\leq J_f(u_{f,0})+\frac{(tR)^2}2 S_2^{\frac N2}-\frac{(tR)^{2^*}}{2^*}S_2^{\frac N2}+\frac{(tR)^2}2(O(\eps^{N-2})+O(\eps^{2-2s}))\\&\quad +\frac{(tR)^{2^*}}{2^*}O(\eps^{  N})-(tR)^{2^*-1} C\eps^{\frac{N-2}2}.
     \end{aligned}
 \end{equation}
For $t=1$, we choose $R_0>0$ such that for any $R\geq R_0$ and $\eps$ small enough,
\begin{equation}\label{(e1)}
    \begin{aligned}
    &\left(\frac{R_0^2}2 S_2^{\frac N2}-\frac{R_0^{2^*}}{2^*}S_2^{\frac N2}-R_0^{2^*-1}C\eps^{\frac{N-2}2}\right) \\
    &\quad +\Bigg(\frac{R_0^2}2(O(\eps^{N-2})+O(\eps^{2-2s}))+\frac{R_0^{2^*}}{2^*}O(\eps^{N})\Bigg)<0.
\end{aligned}
\end{equation}
Thus $J_f(u_{f,0}+RU_\eps)<J_f(u_{f,0})$, for $R\geq R_0$, proving the first estimate of \eqref{energy estimate}. Next, we fix $R=R_0$ and define $$\phi(t):=\frac{(tR_0)^2}2 S_2^{\frac N2}-\frac{(tR_0)^{2^*}}{2^*}S_2^{\frac N2}-(tR_0)^{2^*-1}C\eps^{\frac{N-2}2}.$$ Further, $$N<6-4s \Longleftrightarrow 2-2s>\frac{N-2}{2}.$$ Hence for $t\in[0,1]$, we rewrite \eqref{(d1)} as
\begin{equation}
    J_f(u_{f,0}+tR_0U_{\eps})\leq J_f(u_{f,0})+\phi(t)+o(\eps^{  \frac{N-2}{2}}).
\end{equation}
Notice that $\phi(0)=0$, $\phi(t)>0$ for small $t>0$ and by \eqref{(e1)}, $\phi(1)<0$. So $\phi$ attains its maximum at some $t_\eps\in(0,1)$. Then,
\begin{equation}\label{(f1)}
    J_f(u_{f,0}+tR_0U_{\eps})\leq J_f(u_{f,0})+\phi(t_\eps)+o(\eps^{ \frac{N-2}{2}}).
\end{equation}
Suppose $t_\eps\to0$ as $\eps\to0$. Then for small enough $\eps$, from \eqref{(f1)} we get
$$J_f(u_{f,0}+tR_0U_{\eps})<J_f(u_{f,0})+\frac1N  S_2^{\frac{N}{2}}.$$
If $t_{\eps} \not \rightarrow 0$, then there exists some $T\in(0,1)$ such that $T<t_\eps<1$ for every $\eps$ small. Moreover, we observe that $$\max_{t\geq0}\left(\frac{(tR_0)^2}2 S_2^{\frac N2}-\frac{(tR_0)^{2^*}}{2^*}S_2^{\frac N2}\right)=\frac{(t_0R_0)^2}2 S_2^{\frac N2}-\frac{(t_0R_0)^{2^*}}{2^*}S_2^{\frac N2}=\frac1N S_2^{\frac{N}{2}},$$ where $t_0^{2^*-2}=\frac{R_0^2}{R_0^{2^*}}$. Thus \eqref{(f1)} gives
\begin{equation}\label{fin_est}
    J_f(u_{f,0}+tR_0U_{\eps})\leq J_f(u_{f,0})+\frac1N S_2^{\frac{N}{2}}-(TR_0)^{2^*-1}C\eps^{\frac{N-2}2}+o(\eps^{ \frac{N-2}{2}}).
\end{equation}
Hence, there exists $\eps_0\in(0,1)$ such that for all $\eps<\eps_0$ and any $t\in[0,1]$, $$J_f(u_{f,0}+tR_0U_{\eps})< J_f(u_{f,0})+\frac1NS_2^{\frac{N}{2}}.$$
For sufficiently small $\eps>0$, taking $\Psi=U_\eps\in \WW_2$ we obtain \eqref{energy estimate}. This completes the proof.
\end{proof}
Next, we prove the $\text{(PS)}_c$ condition for the energy functional $J_f$ under the level mentioned in \eqref{energy estimate}.

\begin{lemma}\label{PS cond} Let $f\gneqq0$. Then $J_f$ satisfies $\text{(PS)}_c$ condition for every
    \begin{equation}
        c<J_f(u_{f,0})+\frac1NS_2^{\frac{N}2}.
    \end{equation}
\end{lemma}
\begin{proof}
    Let $\{u_n\}\subset \WW_2$ be a (PS) sequence of $J_f$ at the level $c$, i.e.
    $J_f(u_n)\to c$ and $\|J_f'(u_n)\|_{\WW_2^\ast}\to0 \text{ as }n\to\infty.$ Observe that
    \begin{align*}
        c+C_1\rho_2(u_n)+o_n(1) & =J_f(u_n)-\frac{1}{2^*}J_f'(u_n)(u_n) \\
        & \geq \left(\frac12-\frac1{2^*}\right)\rho_2(u_n)^2-\left(1-\frac1{2^*}\right)\|f\|_{\WW_2^*}\rho_2(u_n).
    \end{align*}
    Hence $\{u_n\}$ is bounded in $\WW_2$. By the reflexivity of $\WW_2$, up to a subsequence, $u_n\rightharpoonup u_0$ in $\WW_2$. Then
    \begin{align}\label{aba1}
        \|\nabla u_n\|_2^2-\|\nabla (u_n-u_0)\|_2^2=\|\nabla u_0\|_2^2+o_n(1).
    \end{align}
    From the compact embedding of $\WW_2 \hookrightarrow L_{\text{loc}}^2(\RR)$ (Proposition \ref{compact embeddings_1}), we see that $u_n(x) \ra u_0(x)$ for a.e. $x \in \RR$. By the Brezis-Lieb lemma,
    \begin{equation}\label{aba}
        \begin{aligned}
        &[u_n]_{s,2}^2-[u_n-u_0]_{s,2}^2=[u_0]_{s,2}^2+o_n(1), \\
        &\|(u_n)_+\|_{2^*}^{2^*}-\|(u_n-u_0)_+\|_{2^*}^{2^*}=\|(u_0)_+\|_{2^*}^{2^*}+o_n(1).
    \end{aligned}
    \end{equation}
    Suppose $\{ u_n\}$ does not converge to $u_0$ i.e. $\rho_2(u_n-u_0)\geq C$ for all $n$. Since $f\in \WW_2^*$, using \eqref{aba1} and \eqref{aba},
    \begin{equation}\label{aba2}
        \rho_2(u_n-u_0)^2-\|(u_n-u_0)_+\|_{2^*}^{2^*}=J_f'(u_n)(u_n-u_0)+o_n(1)=o_n(1).
    \end{equation}
    By the Sobolev inequality $S_2\norm{u_+}^2_{2^*} \le \rho_2(u)^2$ and using \eqref{aba2}, we have
    \begin{align*}
        \|(u_n-u_0)_+\|_{2^*}^{2^*-2}=\frac{\|(u_n-u_0)_+\|_{2^*}^{2^*}}{\|(u_n-u_0)_+\|_{2^*}^2}=\frac{\rho_2(u_n-u_0)^2+o_n(1)}{\|(u_n-u_0)_+\|_{2^*}^2} \geq S_2+\frac{o_n(1)}{\rho_2(u_n-u_0)^2}.
    \end{align*}
    Since $\rho_2(u_n-u_0)\geq C$, $$\|(u_n-u_0)_+\|_{2^*}^{2^*-2}\geq S_2+o_n(1) \Longrightarrow \|(u_n-u_0)_+\|_{2^*}^{2^*}\geq S_2^{\frac{N}2}+o_n(1). $$
    Therefore, using \eqref{aba2}, for large $n$,
    \begin{align}
        \frac1NS_2^{\frac{N}2}& \leq \frac1N\|(u_n-u_0)_+\|_{2^*}^{2^*}+o_n(1) \nonumber\\
        &=\frac12\rho_2(u_n-u_0)^2-\frac1{2^*}\|(u_n-u_0)_+\|_{2^*}^{2^*}+o_n(1)\nonumber\\
        &=J_f(u_n)-J_f(u_0)+o_n(1) < J_f(u_{f,0})+\frac1NS_2^{\frac{N}2}-J_f(u_0).\label{mixed2}
    \end{align}
    Next, we claim that $J_f(u_0)\geq J_f(u_{f,0})$. If this claim holds, then we get a contradiction from \eqref{mixed2}. Hence, $u_n\to u_0$ in $\WW_2$ and $J_f$ satisfies (PS) condition at the level $c$. Now we are left to prove the claim. Note that, if $J_f(u_0)\geq0$, then the claim follows as  $c_{f,0}<0$. Now we assume  $J_f(u_0)<0$. For $t\geq0$, we define
    \begin{equation}\label{para-1}
        g(t)=J_f\left(t\frac{u_0}{\rho_2(u_0)}\right)=\frac12t^2-\frac1{2^*}t^{2^*}\frac{\|(u_0)_+\|_{2^*}^{2^*}}{\rho_2(u_0)^{2^*}}-\frac{t}{\rho_2(u_0)}\prescript{}{\WW_2^*}{\langle}f,u_0{\rangle}_{\WW_2}.
    \end{equation}
    We observe that: near zero $g$ is negative and strictly decreasing, and $g(t)\to-\infty$ as $t\to\infty$. Further by Proposition \ref{first_sol}, $J_f\big|_{\pa B_{r_0}}>\delta_0>0.$ Thus, $g$ is positive in a neighbourhood of $r_0$. It is also easy to see that $g'$ has a unique positive critical point. Combining the above facts, we obtain that $g$ has exactly two positive critical points, say $t_1,t_2$, where $0<t_1<r_0\leq t_2$, $t_1$ is the local minimum, and $t_2$ is the global maximum. Since $g(\rho_2(u_0))=J_f(u_0)<0$, $u_0$ is nonzero. Further since $u_0$ is a critical point of $J_f$, $$g'(\rho_2(u_0))=\left\langle J_f'(u_0),\frac{u_0}{\rho_2(u_0)}\right\rangle=0.$$ Thus $\rho_2(u_0)$ is a positive critical point of $g$. Hence there are two possible scenarios $\rho_2(u_0)=t_1$ or $\rho_2(u_0)=t_2$. Since $g(t_2)>0$, $\rho_2(u_0)$ can not be $t_2$. Thus $\rho_2(u_0)=t_1<r_0$. By Proposition \ref{first_sol}, we infer that $J_f(u_0)\geq c_{f,0}=J_f(u_{f,0})$. This completes the proof of the lemma.
\end{proof}
Finally, we prove the existence of a second positive solution of \eqref{main_prob} using the Mountain pass theorem.
\begin{proof}[{\bf Proof of Theorem \ref{Main_Theorem}:}]
From Proposition \ref{first_sol}, we obtain the first solution $$u_{f,0}\in B_{r_0},\quad J_f\big|_{\pa B_{r_0}}\geq \delta_0>0.$$ By Proposition \ref{energy_prop}, There exists $R_0>0$ (enlarging $R_0$ if necessary) such that
\begin{align*}
&\rho_2(u_{f,0}+R_0U_\eps)>r_0, \; J_f(u_{f,0}+R_0U_\eps)<J_f(u_{f,0})=c_{f,0}<0, \text{ and }\\
    &J_f(u_{f,0}+tR_0U_\eps)<J_f(u_{f,0})+\frac1N S_2^{\frac{N}{2}},\quad\forall\,t\in[0,1].
\end{align*}
Define
\begin{align*}
 & c:=\inf_{\gamma\in\Gamma}\max_{t\in[0,1]}J_f(\gamma(t)), \text{ where } \\ & \Gamma:=\left\{\gamma\in \mathcal{C}([0,1],\WW_2):\gamma(0)=u_{f,0},\,\gamma(1)=u_{f,0}+R_0U_\eps\right\}.
\end{align*}
Observe that the path $\gamma_0(t):=u_{f,0}+tR_0U_\eps,\, t\in[0,1]$ lies in $\Gamma$. Thus, $$c\leq\max_{t\in[0,1]}J_f(\gamma_0(t))<J_f(u_{f,0})+\frac1N S_2^{\frac{N}2}.$$
Therefore, by Lemma \ref{PS cond} and the Mountain pass theorem, there exists a critical point $u_{f,1}$ of $J_f$ such that $J_f(u_{f,1})=c$. Moreover, since $c\geq\delta_0>0>c_{f,0}$, we get $u_{f,1} \neq u_{f,0}$.
\end{proof}

\section{Existence of a positive solution for nonlinear non-homogeneous problem}\label{p_prob_section}

Theorem \ref{non-existence theorem} naturally raises the question of whether one can obtain the existence of a nontrivial solution of \eqref{main_prob_p} under a nontrivial perturbation $f \in \WW_p^*$. The aim of this section is to provide an affirmative answer. More precisely, we establish that \eqref{main_prob_p} admits a positive solution provided $\norm{f}_{\WW_p^*}$ is sufficiently small. We observe that the proof of Proposition \ref{first_sol} is not directly adaptable to the nonlinear setting, since there may not exist $r_0>0$ such that $J_f$ is strictly convex in $B_{r_0}$. Moreover, we have pointed out that even for $p=2$, the (PS) decomposition of the energy functional $J_f$ is unknown. This motivates us to study a concentration-compactness principle for $\del_p +\frap$, due to the non-compactness of $\WW_p \hookrightarrow L^{p^*}(\RR)$. The concentration-compactness principle was introduced by P.L. Lions in \cite[Lemma I.1]{Li}, who initially developed it for bounded sequences in $\mathcal{D}^{1,p}(\RR)$.
\vspace{0.2cm}

\noi \textbf{Setup:} Let $\{u_n\}$ be a bounded sequence in $\WW_p$. Then, up to a subsequence (still denoting by $u_n$) $u_n\rightharpoonup u_0$ in $\WW_p$. Define $$\mu_n=(|\nabla u_n|^p+|D^su_n|^p)\dx,\quad\nu_n=|u_n|^{p^*}\dx.$$ Observe that, $\{\mu_n\}$ and $\{\nu_n\}$ are uniformly bounded in the space of finite Radon measures $\mathcal{M}(\R^N)$. Further, from the Riesz Representation theorem, $\mathcal{M}(\R^N)$ is the topological dual of $\mathcal{C}_0(\RR)$, where $\mathcal{C}_0(\RR)$ is the space of continuous functions that vanish at infinity. By Banach-Alaoglu theorem, we get two non-negative measures $\mu_0$ and $\nu_0$ so that $\mu_n\overset{\ast}{\rightharpoonup}\mu_0$ and $\nu_n\overset{\ast}{\rightharpoonup}\nu_0$ (as the dual of $\mathcal{C}_0(\RR)$) with respect to the weak star topology of $\mathcal{M}(\RR)$.

\begin{proposition}[Concentration-compactness]\label{CC}
Let $\{ u_n \}$ be a bounded sequence in $\WW_p$. Then there exists at most countable index set $J$ and families $$\{x_j\in\RR:j\in J\},\quad \{\mu_j>0:j\in J\},\; \text{ and }\; \{\nu_j>0:j\in J\},$$ such that
\begin{align}
    &\mu_0 \ge (|\nabla u_0|^p+|D^su_0|^p)\dx+\sum_{j\in J}\mu_j\delta_{x_j}, \; \mu_j = \mu_0(\{x_j\}),\label{cc-1}\\
    &\nu_0 = |u_0|^{p^*}\dx+\sum_{j\in J}\nu_j\delta_{x_j}, \; \nu_j = \nu_0(\{x_j\}),\label{cc-2}\\
    &S_p\nu_j^{\frac{p}{p^*}}\leq\mu_j\label{cc-3}.
\end{align}
Moreover,
\begin{align}
    &\limsup_{n\to\infty}\int_{\RR}(|\nabla u_n|^p+|D^su_n|^p)\dx=\int_{\RR}\d\mu_0+\mu_\infty,\label{cc-6}\\
    &\limsup_{n\to\infty}\int_{\RR}|u_n|^{p^*}\dx=\int_{\RR}\d\nu_0+\nu_\infty,\label{cc-7}\\
    & S_p\nu_{\infty}^{\frac{p}{p^*}}\leq\mu_{\infty},\label{cc-8}
\end{align}
where
\begin{align}
    &\mu_{\infty} := \lim_{R \ra \infty} \limsup_{n\to\infty} \int_{\abs{x} > R} (|\nabla u_n|^p+|D^su_n|^p)\dx,\label{cc-4}\\
    &\nu_{\infty} := \lim_{R \ra \infty} \limsup_{n\to\infty} \int_{\abs{x} > R} |u_n|^{p^*}\dx.\label{cc-5}
\end{align}
\end{proposition}

\begin{proof}
\noi \textbf{Step 1:}    Consider $v_n:=u_n-u_0\rightharpoonup0$ in $\WW_p$. Using the Brezis-Lieb lemma, we write $$\om_n:=\nu_n-|u_0|^{p^*}\dx=(|u_n|^{p^*}-|u_0|^{p^*})\dx=|u_n-u_0|^{p^*}\dx+o_n(1)=|v_n|^{p^*}\dx+o_n(1).$$
    Suppose $\la_n:=(|\nabla v_n|^p+|D^sv_n|^p)\dx \overset{\ast}{\rightharpoonup} \la_0$ and $\om_n \overset{\ast}{\rightharpoonup} \om_0$. Clearly $\la_0,\om_0\geq0$ and $\om_0=\nu_0-|u_0|^{p^*}\dx$. Now for any $\phi\in \mathcal{C}_c^{\infty}(\RR)$, by Sobolev inequality, we see that
    \begin{align*}
        &S_p\left(\int_{\RR}|\phi|^{p^*}\d\om_0\right)^{\frac p{p^*}}=S_p\lim_{n\to\infty}\left(\int_{\RR}|v_n\phi|^{p^*}\dx\right)^{\frac p{p^*}}\\
        &\leq \lim_{n\to\infty}\int_{\RR}(|\nabla(v_n\phi)|^{p}+|D^s(v_n\phi)|^p)\dx\\
        &\leq \lim_{n\to\infty}\int_{\RR}|\phi|^p(|\nabla v_n|^{p}+(1+\theta)|D^sv_n|^p)\dx \\
        & \quad +\int_{\RR}|v_n|^p(|\nabla \phi|^{p}+C_\theta|D^s\phi|^p)\dx\\& \quad +C\int_{\RR}\phi^{p-1}v_n|\nabla v_n|^{p-1}|\nabla\phi|+ v_n^{p-1}\phi|\nabla\phi|^{p-1}|\nabla v_n|\dx\\
        &\leq \lim_{n\to\infty}\int_{\RR}|\phi|^p(|\nabla v_n|^{p}+(1+\theta)|D^sv_n|^p)\dx \\ & \quad +\int_{\RR}|v_n|^p(|\nabla \phi|^{p}+C_\theta|D^s\phi|^p)\dx\\ & \quad +C(\|\phi \abs{\nabla v_n}\|_p^{p-1}\|v_n|\nabla\phi|\|_p+\|v_n|\nabla \phi|\|_p^{p-1}\|\phi|\nabla v_n|\|_p)\\
        &\leq \lim_{n\to\infty}\int_{\RR}|\phi|^p|\nabla v_n|^{p}\dx+(1+\theta)\lim_{n\to\infty}\int_{\RR}|\phi|^p|D^sv_n|^p\dx,\; \forall \, \theta>0.
    \end{align*}
    Taking $\theta \ra 0$,
\begin{align*}
    S_p\left(\int_{\RR}|\phi|^{p^*}\d\om_0\right)^{\frac p{p^*}} \le \lim_{n\to\infty}\int_{\RR}|\phi|^p\d\la_n=\int_{\RR}|\phi|^p\d\la_0.
\end{align*}
The second inequality comes by using \eqref{Leibniz} and the inequality
\begin{align}\label{ele}
    (1+a)^t\leq 1+a^t+C(a+a^{t-1}),a>0,t\geq1.
\end{align}
    To get the last inequality, we use the facts that $\{v_n\}$ is bounded in $\WW_p$, $v_n\rightharpoonup0$ in $\WW_p$, and the compact embeddings of $\WW_p\hookrightarrow L^p_{\text{loc}}(\RR)$ and $\WW_p\hookrightarrow L^p(|D^s\phi|^p,\RR).$ Hence
    \begin{equation}
        S_p^{\frac{1}{p}}\left(\int_{\RR}|\phi|^{p^*}\d\om_0\right)^{\frac{1}{p^*}}\leq \left(\int_{\RR}|\phi|^p\d\la_0\right)^{\frac{1}p}, \; \forall \, \phi\in \mathcal{C}_c^{\infty}(\RR).
    \end{equation}
    Therefore by \cite[Lemma 1.2]{Li}, there exist an at most countable set $J$ and families $\{x_j\}_{j\in J}$ in $\RR$ and $\{\nu_j\}_{j\in J}$ in $(0,\infty)$ such that
    \begin{equation}
        \om_0=\sum_{j\in J}\nu_j\de_{x_j}.
    \end{equation}
    This implies that \eqref{cc-2} holds. From the weak convergence $u_n\rightharpoonup u_0$ and the weak lower semicontinuity of $\rho_p(\cdot)$, we have $\mu_0\geq(|\nabla u_0|^p+|D^s u_0|^p)\dx$. Since $(|\nabla u_0|^p+|D^s u_0|^p)\dx$ and $\sum_{j\in J}\mu_j\de_{x_j}$ are mutually singular, we obtain
    $$\mu_0\geq \left(|\nabla u_0|^p+|D^s u_0|^p\right)\dx+\sum_{j\in J}\mu_j\de_{x_j}.$$

 \noi \noi \textbf{Step 2:}   Let $0\leq \phi\leq1, \phi(0)=1, \supp{\phi} \subset B_1(0)$ and define $\phi_{j,\eps}(x)=\phi\left(\frac{x-x_j}{\eps}\right)$. We also have
    \begin{align*}
        &S_p\left(\int_{\RR}|\phi_{j,\eps}|^{p^*}\d\nu_0\right)^{\frac p{p^*}}=S_p\lim_{n\to\infty}\left(\int_{\RR}|u_n\phi_{j,\eps}|^{p^*}\dx\right)^{\frac p{p^*}}\\
        &\leq \lim_{n\to\infty}\left(\int_{\RR}(|\nabla(u_n\phi_{j,\eps})|^{p}+|D^s(u_n\phi_{j,\eps})|^p\dx\right)\\
        &\le\lim_{n\to\infty}\int_{\RR}|\phi_{j,\eps}|^p(|\nabla u_n|^{p}+(1+\theta)|D^su_n|^p)\dx \\ & \quad +\int_{\RR}|u_0|^p(|\nabla \phi_{j,\eps}|^{p}+C_\theta|D^s\phi_{j,\eps}|^p)\dx \\ & \quad +C\eps^{-1}\left(\int_{\RR}|u_0|^p\left|\nabla \phi\left(\frac{x-x_j}{\eps}\right)\right|^p\dx\right)^{\frac1p} \\ & \quad +C\eps^{1-p} \left(\int_{\RR}|u_0|^p\left|\nabla \phi\left(\frac{x-x_j}{\eps}\right)\right|^p\dx\right)^{\frac{p-1}p},
        \end{align*}
        where we use \eqref{Leibniz}, \eqref{ele}, and the compactness of $\mathcal{W}_p \hookrightarrow L_{\text{loc}}^p(\RN)$ for the second inequality. Now,  using the change of variables, we estimate
        \begin{align*}
            &\int_{\RR}|u_0|^p|\nabla \phi_{j,\eps}|^{p}\dx=\eps^{-p}\int_{B_{\eps}(x_j)}|u_0|^p\left|\nabla \phi\left(\frac{x-x_j}{\eps}\right)\right|^p\dx\\&\leq\eps^{-p}\left(\int_{B_{\eps}(x_j)}|u_0|^{p^*}\right)^{\frac{p}{p^*}}\left(\int_{B_{\eps}(x_j)}\left|\nabla \phi\left(\frac{x-x_j}{\eps}\right)\right|^{\frac{pp^*}{p^*-p}}\dx\right)^{\frac{p^*-p}{p^*}}\\
            &=\left(\int_{B_{\eps}(x_j)}|u_0|^{p^*}\right)^{\frac{p}{p^*}}\left(\int_{B_1(0)}\left|\nabla \phi(x)\right|^{\frac{pp^*}{p^*-p}}\dx\right)^{\frac{p^*-p}{p^*}}\\
            &\leq C\left(\int_{B_{\eps}(x_j)}|u_0|^{p^*}\right)^{\frac{p}{p^*}}\to0, \text{ as } \eps \ra 0,
        \end{align*}
       and
        \begin{align*}
            &\eps^{-1}\left(\int_{\RR}|u_0|^p\left|\nabla \phi\left(\frac{x-x_j}{\eps}\right)\right|^p\dx\right)^{\frac1p}+\eps^{1-p} \left(\int_{\RR}|u_0|^p\left|\nabla \phi\left(\frac{x-x_j}{\eps}\right)\right|^p\dx\right)^{\frac{p-1}p}\\
            &\leq C\left(\int_{B_{\eps}(x_j)}|u_0|^{p^*}\right)^{\frac{1}{p^*}}+\left(\int_{B_{\eps}(x_j)}|u_0|^{p^*}\right)^{\frac{p-1}{p^*}}\to0, \text{ as } \eps \ra 0.
        \end{align*}
        Further, using the arguments for \cite[Eq. (2.4)]{BoSaSi}, the following convergences holds:
        \begin{align*}
            \int_{\RR}|u_0|^p|D^s\phi_{j,\eps}|^{p}\dx\to0, \text{ as } \eps \ra 0.
        \end{align*}
        Hence
        \begin{align*}
            S_p\left(\int_{\RR}|\phi_{j,\eps}|^{p^*}\d\nu_0\right)^{\frac p{p^*}} \leq \lim_{n\to\infty}\int_{\RR}|\phi_{j,\eps}|^p(|\nabla u_n|^{p}+(1+\theta)|D^su_n|^p)\dx+o_\eps(1),\; \forall \,\theta>0.
        \end{align*}
        Taking $\theta \ra 0$,
        \begin{align}\label{ele-1}
            S_p\left(\int_{\RR}|\phi_{j,\eps}|^{p^*}\d\nu_0\right)^{\frac p{p^*}} \le \lim_{n\to\infty}\int_{\RR}|\phi_{j,\eps}|^p\d\mu_n+o_\eps(1) =\int_{\RR}|\phi_{j,\eps}|^p\d\mu_0+o_\eps(1).
        \end{align}
        Now,
    $$\int_{\RR}|\phi_{j,\eps}|^{p^*}\d\nu_0\geq\int_{\RR}|\phi_{j,\eps}|^{p^*}\nu_j \, {\rm d}\de_{x_j} = \nu_j, \, \text{ for every } \, \eps >0,$$ and $$\int_{\RR}|\phi_{j,\eps}|^{p^*}\d\mu_0\leq\mu_0(B_\eps(x_j))\to\mu_0(\{x_j\}), \; \text{as } \; \eps\to0.$$
        Taking $\eps\to0$ in \eqref{ele-1}, we get $$S_p\nu_j^{\frac p{p^*}}\leq\mu_j, \; \forall \, j\in J.$$
     Thus, \eqref{cc-3} holds.

  \noi \textbf{Step 3:}  This step is inspired by \cite[Proposition 2]{Chabrowski1995}, providing the information about a possible loss of mass at infinity.  For $R>0$, set $\psi_R \in \mathcal{C}_b(\rd)$, defined as $\psi_R(x) = 1$ if $\abs{x}>R$, $\psi_R(x) = 0$ if $\abs{x} \le \frac{R}{2}$, and $|\nabla\psi_R|\leq\frac CR$. We write
\begin{align}\label{cc8}
    &\limsup_{n \ra \infty} \int_{\rd} \left( |\nabla u_n|^p + \abs{D^s u_n}^p \right)\dx \no \\
    &= \limsup_{n \ra \infty} \left( \int_{\rd} \psi_R^p \left( |\nabla u_n|^p + \abs{D^s u_n}^p \right) + \left( 1-\psi_R^p \right) \left( |\nabla u_n|^p + \abs{D^s u_n}^p\right)\dx\right),
\end{align}
   Since $(1-\psi_R^p)\in \mathcal{C}_0(\RR)$, $$\limsup_{n \ra \infty} \int_{\rd} \left( 1-\psi_R^p \right) \left( |\nabla u_n|^p + \abs{D^s u_n}^p \right) \dx=\int_{\RR}(1-\psi_R^p)\d\mu_0.$$
   Notice that $\psi_R\to0$ pointwise and $\mu_0$ is a nonnegative measure. By dominated convergence theorem, $\lim_{R\to\infty}\int_{\RR}\psi_R^p\d\mu_0=0.$ Hence, we get
   \begin{equation}\label{cc9}
       \lim_{R\to\infty}\limsup_{n \ra \infty} \int_{\rd} \left( 1-\psi_R^p \right) \left( |\nabla u_n|^p + \abs{D^s u_n}^p \right) \dx=\int_{\RR}\d\mu_0.
   \end{equation}
   Observe that
   \begin{align*}
       \int_{\{|x|>R\}} \left( |\nabla u_n|^p + \abs{D^s u_n}^p \right) \dx&\leq \int_{\RR}\psi_R^p \left( |\nabla u_n|^p + \abs{D^s u_n}^p \right) \dx\\&\leq \int_{\{|x|\geq R/2\}} \left( |\nabla u_n|^p + \abs{D^s u_n}^p \right) \dx.
   \end{align*}
   By Sandwich theorem,
   \begin{align}\label{cc10}
      \mu_{\infty}&=\lim_{R\to\infty}\limsup_{n\to\infty}\int_{\{|x|>R\}}\left( |\nabla u_n|^p + \abs{D^s u_n}^p \right) \dx\no\\&=\lim_{R\to\infty}\limsup_{n\to\infty} \int_{\rd} \psi_R^p \left( |\nabla u_n|^p + \abs{D^s u_n}^p \right) \dx.
   \end{align}
   Combining \eqref{cc8}-\eqref{cc10}, we get \eqref{cc-6}. Similarly we get \eqref{cc-7} and $$\nu_\infty=\lim_{R\to\infty}\limsup_{n\to\infty}\int_{\RR}|u_n\psi_R|^{p^*}\dx.$$By the definition of $S_p$ and using \eqref{Leibniz}, \eqref{ele} for any $\theta>0$, we get
   \begin{align*}
       &S_p\left(\int_{\RR}|u_n\psi_R|^{p^*}\dx\right)^{\frac{p}{p^*}}\leq \int_{\RR}|\nabla (u_n\psi_R)|^p+|D^s(u_n\psi_R)|^p\dx\\
       &\leq \int_{\RR}|\psi_R|^p(|\nabla u_n|^{p}+(1+\theta)|D^su_n|^p)\dx+\int_{\RR}|u_n|^p(|\nabla \psi_R|^{p}+C_\theta|D^s\psi_R|^p)\dx\\&\quad+C\left(\int_{\RR}|u_n|^p\left|\nabla \psi_R\right|^p\dx\right)^{\frac1p}+ C\left(\int_{\RR}|u_n|^p\left|\nabla \psi_R\right|^p\dx\right)^{\frac{p-1}p},
   \end{align*}
   Using the fact that $\nabla\psi_R$ has compact support, $\WW_p\Subset L^p_{\text{loc}}$ and $|\nabla\psi_R|\leq\frac CR$, we see that
   \begin{align*}
    \lim_{R\to\infty}\limsup_{n\to\infty}\int_{\RR}|u_n|^p|\nabla \psi_R|^{p}\dx 
    +C\left(\int_{\RR}|u_n|^p\left|\nabla \psi_R\right|^p\dx\right)^{\frac1p} 
    + C\left(\int_{\RR}|u_n|^p\left|\nabla \psi_R\right|^p\dx\right)^{\frac{p-1}p}=0.
   \end{align*}
   Following the arguments of \cite[Eq. (2.9)]{BoSaSi}, $$\lim_{R\to\infty}\limsup_{n\to\infty}\int_{\RR}|u_n|^p|D^s\psi_R|^p\dx=0.$$
   Hence, for any $\theta>0$,
   $$S_p\nu_{\infty}^{\frac{p}{p^*}}\leq \lim_{R\to\infty}\limsup_{n\to\infty}\int_{\RR}\psi_R^p(|\nabla u_n|^{p}+(1+\theta)|D^su_n|^p)\dx.$$ Finally taking $\theta\to0$, we obtain \eqref{cc-8}.
\end{proof}
As an application of the concentration-compactness principle, we prove the following.

\begin{proposition}\label{grad_p.w.}
   Let $\{ u_n \}$ be a bounded sequence in $\WW_p$ such that $J_f(u_n)\to c$ and $J_f'(u_n)\to0$ in $\WW_p^*$. Assume that $u_n\rightharpoonup u_0$ in $\WW_p$. Then the set $J$ in Proposition \ref{CC} is finite. Further, the following convergences hold
   \begin{align*}
    &\frac{\pa u_n}{\pa x_k}\to \frac{\pa u_0}{\pa x_k}\text{ a.e. in }\RR,\\
    &|\nabla u_n|^{p-2}\nabla u_n\rightharpoonup |\nabla u_0|^{p-2}\nabla u_0 \text{ in }([L^{p}(\RR)]^N)^*.
\end{align*}
\end{proposition}
\begin{proof}
We follow the approach of \cite[Lemma 2.2]{Ji}. For $\eps>0$, $\phi \in \mathcal{C}_c^{\infty}(\RR)$ with $0\leq \phi\leq1,\phi(0)=1, \supp{\phi} \subset B_1(0)$, we consider the function $\phi_{j,\eps}(x)=\phi\left(\frac{x-x_j}{\eps}\right)$. Since $\langle J_f'(u_n),\phi\rangle\to0$ for every $\phi\in\WW_p$, taking $\phi_{j,\eps}u_n$ as a test function,
\begin{equation}\label{mixed1}
\begin{aligned}
    &\int_{\RR}|\nabla u_n|^{p-2}\nabla u_n\cdot\nabla(\phi_{j,\eps}u_n)\dx+\AA_p(u_n,\phi_{j,\eps}u_n)\\
    &=\int_{\RR}|u_n|^{p^*}\phi_{j,\eps}\dx+{}_{\WW_p^*}\langle f,u_n\phi_{j,\eps} \rangle_{\WW_p}+o_n(1).
\end{aligned}
\end{equation}
By the compact embedding $\WW_p\hookrightarrow L^p_{\text{loc}}(\RR)$, $u_n(x)\rightarrow u_0(x)$ a.e. $x \in \RR$. Since $\{u_n\phi_{j,\eps}\}$ is bounded in $\WW_p$, we infer that $u_n\phi_{j,\eps}\rightharpoonup u_0\phi_{j,\eps}$ in $\WW_p$. Thus ${}_{\WW_p^*}\langle f,u_n\phi_{j,\eps} \rangle_{\WW_p}\to{}_{\WW_p^*}\langle f,u_0\phi_{j,\eps} \rangle_{\WW_p}.$ Since  $u_n\rightharpoonup u_0$ in $\WW_p$, we get $\nu_n\overset{\ast}{\rightharpoonup}\nu_0$, and hence  $$\int_{\RR}|u_n|^{p^*}\phi_{j,\eps}\dx=\int_{\RR}\phi_{j,\eps}\d\nu_n\to\int_{\RR}\phi_{j,\eps}\d\nu_0, \; \text{ as } n \ra \infty.$$
Expanding the left hand side of \eqref{mixed1}, we get
\begin{align*}
    &\int_{\RR}|\nabla u_n|^{p-2}\nabla u_n\cdot\nabla(\phi_{j,\eps}u_n)\dx+\AA_p(u_n,\phi_{j,\eps}u_n)\\
    &=\int_{\RR}u_n|\nabla u_n|^{p-2}\nabla u_n\cdot\nabla\phi_{j,\eps}\dx+\int_{\RR}\phi_{j,\eps}|\nabla u_n|^p\dx\\
    &\quad+\iint_{\R^{2N}}\frac{A_p(u_n)(x,y)(\phi_{j,\eps}(x)-\phi_{j,\eps}(y))u_n(x)}{|x-y|^{N+sp}}\dxy+\int_{\RR}\phi_{j,\eps}|D^su_n|^p\dx\\
    &=\int_{\RR}\phi_{j,\eps}\d\mu_n+\int_{\RR}u_n|\nabla u_n|^{p-2}\nabla u_n\cdot\nabla\phi_{j,\eps}\dx\\
    &\quad+\iint_{\R^{2N}}\frac{A_p(u_n)(x,y)(\phi_{j,\eps}(x)-\phi_{j,\eps}(y))u_n(x)}{|x-y|^{N+sp}}\dx\dy.
\end{align*}
Moreover, for every $\eps$, applying H\"{o}lder's inequality and the compact embeddings $\WW_p\hookrightarrow L^p_{\text{loc}}(\RR),\WW_p\hookrightarrow L^p(|D^s\phi_{j,\eps}|^p;\RR)$, we obtain
\begin{align*}
    &\lim_{n\to\infty}\left|\int_{\RR}u_n|\nabla u_n|^{p-2}\nabla u_n\cdot\nabla\phi_{j,\eps}+\iint_{\R^{2N}}\frac{A_p(u_n)(x,y)(\phi_{j,\eps}(x)-\phi_{j,\eps}(y))u_n(x)}{|x-y|^{N+sp}}\right|\\
    &\leq\lim_{n\to\infty} \left( \|\nabla u_n\|_p^{p-1}\|u_n|\nabla\phi_{j,\eps}|\|_p+[u_n]_{s,p}^{p-1}\|u_n|D^s\phi_{j,\eps}|\|_p \right)\\
    &\leq C \left(\|u_0|\nabla\phi_{j,\eps}|\|_p + \|u_0|D^s\phi_{j,\eps}|\|_p\right).
\end{align*}
In view of the proof of Proposition \ref{CC}, we similarly get $$\|u_0|\nabla\phi_{j,\eps}|\|_p+\|u_0|D^s\phi_{j,\eps}|\|_p=o_\eps(1).$$
Hence
\begin{align}\label{eq_eps1.1}
    {}_{\WW_p^*}\langle f,u_0\phi_{j,\eps} \rangle_{\WW_p}+\int_{\RR}\phi_{j,\eps}\d\nu_0 =\lim_{n\to\infty}\int_{\RR}\phi_{j,\eps}\d\mu_n + o_{\eps}(1) =\int_{\RR}\phi_{j,\eps}\d\mu_0 + o_{\eps}(1).
\end{align}
Again, using a similar set of arguments as in Proposition \ref{CC},
$$\left|{}_{\WW_p^*}\langle f,u_0\phi_{j,\eps} \rangle_{\WW_p}\right|\leq \|f\|_{\WW_p^*}\rho_p(u_0\phi_{j,\eps})=o_\eps(1).$$
Taking $\eps\to0$ in \eqref{eq_eps1.1}, we obtain
$$\nu_j=\mu_j.$$
Since $S_p\nu_j^{\frac{p}{p^*}}\leq\mu_j$, we have
\begin{equation}\label{mixed3}
    \mu_j=\nu_j\geq S_p^{\frac{N}{p}}, \text{ where } \sum_j\nu_j<\infty.
\end{equation}
Now \eqref{mixed3} infers that the set $J$ is finite. Suppose $J=\{1,\ldots,m\}$. Choose $\eps,R>0$ be such that $B_\eps(x_i)\cap B_\eps(x_j)=\emptyset$ for every $1\leq i\neq j\leq m$ and $\cup_{j=1}^m B_{\eps}(x_j)\subset B_R(0)$. Let $\xi_\eps\in \mathcal{C}_c^{\infty}(B_{2R}(0))$ be such that $0\le \xi_\eps \le 1$, $\xi_\eps\equiv0$ on $\cup_{j=1}^m B_{\frac12\eps}(x_j)$ and $\xi_\eps \equiv 1$ on $B_R(0) \setminus \cup_{j=1}^m B_{\eps}(x_j)$.
Since $\langle J_f'(u_n),\xi_\eps(u_n-u_0)\rangle\to0$, we see that
\begin{align*}
    o_n(1)&=\int_{\RR}|\nabla u_n|^{p-2}\nabla u_n\cdot\nabla(\xi_\eps(u_n-u_0))\dx+\AA_p(u_n,\xi_\eps(u_n-u_0))\\
    &\quad-\int_{\RR}|u_n|^{p^*-2}u_n(\xi_\eps(u_n-u_0))\dx-{}_{\WW_p^*}\langle f,\xi_\eps(u_n-u_0)\rangle_{\WW_p}.
\end{align*}
Notice that the sequence $\{|u_n|^{p^*-2}u_n\}_n$ is bounded in $L^{(p^*)'}(\RR)$ and since $u_n(x)\to u_0(x)$ pointwise a.e. $x \in \RR$, $|u_n|^{p^*-2}u_n\rightharpoonup |u_0|^{p^*-2}u_0$ in $L^{(p^*)'}(\RR)$. Hence $\int_{\RR}|u_n|^{p^*-2}u_n\xi_\eps u_0\dx\to\int_{\RR}|u_0|^{p^*}\xi_\eps\dx$ as $\xi_{\eps}u_0\in L^{p^*}(\RR)$.  Therefore, using Proposition \ref{CC},
\begin{align*}
   \lim_{n \ra \infty} \int_{\RR}|u_n|^{p^*-2}u_n(\xi_\eps(u_n-u_0))\dx&=\lim_{n \ra \infty}\int_{\RR}|u_n|^{p^*}\xi_\eps\dx-\int_{\RR}|u_0|^{p^*}\xi_\eps\dx\\
    &=\lim_{n \ra \infty}\int_{\RR}\xi_\eps(\d\nu_n-\d\nu_0) \\
    &=\sum_{j=1}^m\xi_\eps(x_j)\nu_j=0.
\end{align*}
Further, ${}_{\WW_p^*}\langle f,\xi_\eps u_n\rangle_{\WW_p}\to {}_{\WW_p^*}\langle f,\xi_\eps u_0\rangle_{\WW_p}$. Hence,
$$\int_{\RR}|\nabla u_n|^{p-2}\nabla u_n\cdot\nabla(\xi_\eps(u_n-u_0))\dx+\AA_p(u_n,\xi_\eps(u_n-u_0))=o_n(1).$$
Now
\begin{align*}
    &\int_{\RR}|\nabla u_n|^{p-2}\nabla u_n\cdot\nabla(\xi_\eps(u_n-u_0))\dx\\&=\int_{\RR}(u_n-u_0)|\nabla u_n|^{p-2}\nabla u_n\cdot\nabla\xi_\eps+\xi_{\eps}|\nabla u_n|^{p-2}\nabla u_n\cdot\nabla(u_n-u_0)\dx,
\end{align*}
and
\begin{align*}
    &\AA_p(u_n,\xi_\eps(u_n-u_0))\\
    &=\iint_{\R^{2N}}A_p(u_n)(x,y)(\xi_\eps(u_n-u_0)(x)-\xi_\eps(u_n-u_0)(y))\frac{\dx\dy}{|x-y|^{N+sp}}\\
    &=\iint_{\R^{2N}}A_p(u_n)(x,y)(\xi_\eps(x)-\xi_\eps(y))(u_n(x)-u_0(x))\frac{\dx\dy}{|x-y|^{N+sp}}\\
    &\quad+\iint_{\R^{2N}}\xi_\eps(y)A_p(u_n)(x,y)((u_n-u_0)(x)-(u_n-u_0)(y))\frac{\dx\dy}{|x-y|^{N+sp}}\\
    &=\iint_{\R^{2N}}A_p(u_n)(x,y)(\xi_\eps(x)-\xi_\eps(y))(u_n(x)-u_0(x))\frac{\dx\dy}{|x-y|^{N+sp}}\\
    &\quad+\iint_{\R^{2N}}\xi_\eps(y)|u_n(x)-u_n(y)|^{p}-\xi_\eps(y)A_p(u_n)(x,y)(u_0(x)-u_0(y))\frac{\dx\dy}{|x-y|^{N+sp}}.
\end{align*}
Hence,
\begin{equation}\label{mixed4}
    \begin{aligned}
    &\int_{\RR}\xi_\eps(x)(|\nabla u_n|^p-|\nabla u_n|^{p-2}\nabla u_n\cdot\nabla u_0)\dx\\&+\iint_{\R^{2N}}\frac{\xi_\eps(x)(|u_n(x)-u_n(y)|^p-(A_p(u_n)(x,y)(u_0(x)-u_0(y))))}{|x-y|^{N+sp}}\dx\dy\\
    &=\int_{\RR}(u_0-u_n)|\nabla u_n|^{p-2}\nabla u_n\cdot\nabla\xi_\eps\dx\\&\quad+\iint_{\R^{2N}}\frac{A_p(u_n)(x,y)(\xi_\eps(x)-\xi_\eps(y))(u_0(x)-u_n(x))}{|x-y|^{N+sp}}\dx\dy+o_n(1).
\end{aligned}
\end{equation}
Now using the monotonicity of the function $t\mapsto |t|^{p-2}t$, we get
\begin{align*}
    &0\leq \int_{\RR}\xi_\eps(|\nabla u_n|^{p-2}\nabla u_n-|\nabla u_0|^{p-2}\nabla u_0)\cdot(\nabla u_n-\nabla u_0)\dx\\&\leq
    \int_{\RR}\xi_\eps(|\nabla u_n|^{p-2}\nabla u_n-|\nabla u_0|^{p-2}\nabla u_0)\cdot(\nabla u_n-\nabla u_0)\dx\\&\quad+\iint_{\R^{2N}}\frac{\xi_\eps(x)(A_p(u_n)(x,y)-A_p(u_0)(x,y))((u_n-u_0)(x)-(u_n-u_0)(y))}{|x-y|^{N+sp}}\dx\dy\\
    &=\int_{\RR}\xi_\eps(x)(|\nabla u_n|^p+|\nabla u_0|^p-|\nabla u_n|^{p-2}\nabla u_n\cdot\nabla u_0-|\nabla u_0|^{p-2}\nabla u_0\cdot\nabla u_n)\dx\\
    &\quad+ \iint_{\R^{2N}}\xi_\eps(x)(|u_n(x)-u_n(y)|^p+|u_0(x)-u_0(y)|^p)\frac{\dx\dy}{|x-y|^{N+sp}}\\
    &\quad- \iint_{\R^{2N}}\frac{\xi_\eps(x)(A_p(u_n)(x,y)(u_0(x)-u_0(y))+A_p(u_0)(x,y)(u_n(x)-u_n(y))}{|x-y|^{N+sp}}\dx\dy\\
    &=\int_{\RR}\xi_\eps(x)(|\nabla u_n|^p-|\nabla u_n|^{p-2}\nabla u_n\cdot\nabla u_0)\dx\\
    &\quad+\iint_{\R^{2N}}\frac{\xi_\eps(x)(|u_n(x)-u_n(y)|^p-(A_p(u_n)(x,y)(u_0(x)-u_0(y))))}{|x-y|^{N+sp}}\dx\dy+o_n(1)\\
    &=\int_{\RR}(u_0-u_n)|\nabla u_n|^{p-2}\nabla u_n\cdot\nabla\xi_\eps\dx \\
    &\quad+\iint_{\R^{2N}}A_p(u_n)(x,y)(\xi_\eps(x)-\xi_\eps(y))(u_0(x)-u_n(x))\frac{\dx\dy}{|x-y|^{N+sp}}+o_n(1).
\end{align*}
We use \eqref{mixed4} to get the last equality. Now,
\begin{align*}
    &\left|\int_{\RR}(u_0-u_n)|\nabla u_n|^{p-2}\nabla u_n\cdot\nabla\xi_\eps\dx\right|\leq \|\nabla u_n\|_p^{p-1}\int_{\RR}|u_0-u_n|^p|\nabla\xi_\eps|^p\dx\\
    &\leq C\int_{B_{2R}(0)}|u_0-u_n|^p|\nabla\xi_\eps|^p\dx\to0,\text{ as }n\to\infty,
\end{align*}
using the compactness of $\WW_p\hookrightarrow L^p_{\text{loc}}(\RR)$. Next, we check
\begin{align*}
    &\left|\iint_{\R^{2N}}A_p(u_n)(x,y)(\xi_\eps(x)-\xi_\eps(y))(u_0(x)-u_n(x))\frac{\dx\dy}{|x-y|^{N+sp}}\right|\\
    &\leq [u_n]_{s,p}^{p-1}\int_{\RR}|u_0(x)-u_n(x)|^p|D^s\xi_\eps|^p(x)\dx\\
    &\leq C\int_{\RR}|u_0(x)-u_n(x)|^p|D^s\xi_\eps|^p(x)\dx\to0\text{ as }n\to\infty,
\end{align*}
using the compactness of $\WW_p\hookrightarrow L^p(|D^s\xi_{\eps}|^p,\RR)$.
Hence for any $\eps>0$ small enough, we get
$$0\leq \lim_{n\to\infty}\int_{\RR}\xi_\eps(|\nabla u_n|^{p-2}\nabla u_n-|\nabla u_0|^{p-2}\nabla u_0)\cdot(\nabla u_n-\nabla u_0)\dx\leq0.$$
Let $\eps_0>0$ be arbitrary. For any $\eps<\eps_0$,
\begin{align*}
    0&\leq \lim_{n\to\infty}\int_{B_R(0)\setminus \cup_{j=1}^m B_{\eps_0}(x_j)}(|\nabla u_n|^{p-2}\nabla u_n-|\nabla u_0|^{p-2}\nabla u_0)\cdot(\nabla u_n-\nabla u_0)\dx\\
    &\leq \lim_{n\to\infty}\int_{\RR}\xi_\eps(|\nabla u_n|^{p-2}\nabla u_n-|\nabla u_0|^{p-2}\nabla u_0)\cdot(\nabla u_n-\nabla u_0)\dx \leq0.
\end{align*}
Thus $$(|\nabla u_n|^{p-2}\nabla u_n-|\nabla u_0|^{p-2}\nabla u_0)\cdot(\nabla u_n-\nabla u_0)\to0\text{ a.e. in }B_R(0)\setminus\cup_{j=1}^m B_{\eps_0}(x_j),$$
for every $\eps_0$ small. Hence
$$(|\nabla u_n|^{p-2}\nabla u_n-|\nabla u_0|^{p-2}\nabla u_0)\cdot(\nabla u_n-\nabla u_0)\to0\text{ a.e. in }B_R(0),$$ for any $R>0$ large. From the above convergence, we infer
\begin{align*}
    &\frac{\pa u_n}{\pa x_k}\to \frac{\pa u_0}{\pa x_k}\text{ a.e. in } B_R(0),\; 1\leq k\leq N,\\
    &|\nabla u_n|^{p-2}\frac{\pa u_n}{\pa x_k}\rightharpoonup|\nabla u_0|^{p-2}\frac{\pa u_0}{\pa x_k} \text{ in } (L^p(B_R(0)))^*, \; 1\leq k\leq N.
\end{align*}
To see the first convergence we use the inequalities from \cite[Lemma 3.6]{BaLi}, and the second convergence follows by noticing $\{|\nabla u_n|^{p-2}\frac{\pa u_n}{\pa x_k}\}_n$ is bounded in $(L^p(B_R(0)))^*$.
Finally, taking $R\to\infty$, we get the required result.
\end{proof}

\noi \textbf{Proof of Theorem \ref{Main Theorem_1}:}
By the Sobolev embedding $\WW_p \hookrightarrow L^{p^*}(\RR)$, we can choose $r_0>0$ and $\delta_0>0$  small enough such that
    \begin{equation}\label{lbound-1}
    \begin{aligned}
        &\frac1p\rho_p(u)^p-\frac1{\pst}\|u\|_{\pst}^{\pst}\geq\frac1p\rho_p(u)^p-C_1\rho_p(u)^{\pst}\geq \left\{\begin{array}{ll}
             0, & \forall\,u\in B_{r_0}; \\
             2\delta_0, & \forall\,u\in \pa B_{r_0}.  \\
             \end{array} \right.
    \end{aligned}
    \end{equation}
Thus, $J_f(u)\geq -{}_{\WW_p^*}\langle f,u\rangle_{\WW_p}\geq-\|f\|_{\WW_p^*}r_0>-\infty$ for every $u\in B_{r_0}$. Thus $J_f$ is bounded from below on $B_{r_0}$. Define
\begin{align}\label{min}
    c_{f,1}:=\inf_{B_{r_0}}J_f(u).
\end{align}
For any $v\neq0$ in $\WW_p$ and $t>0$,
    \begin{align*}
        J_f(tv)&<\frac{t^p}p\rho_p(v)^p-t{}_{\WW_p^*}\langle f,v\rangle_{\WW_p}.
    \end{align*}
It is easy to see that, since $f\gneqq 0$, there exists a nonnegative $v\in\WW_p$ such that $J_f(tv)<0$ for sufficiently small $t>0$. Hence $-\infty<c_{f,1}<0$. Now we show that the minimizer of $c_{f,1}$ in \eqref{min} is attained.
Consider a minimizing sequence $u_{n}\in B_{r_0}$ i.e. $J_f(u_{n})\to c_{f,1} \text{ as } n\to\infty.$
We claim that $\rho_p(u_{n})\leq {r_0}-\eps_0$ for every $n$ and for some $\eps_0$ independent of $n$. If not, $\rho_p(u_{n})\to {r_0}$. Thus, by the Sobolev inequality,
\begin{align*}
 c_{f,1}=\lim_{n\to\infty}J_f(u_{n}) & \geq \lim_{n\to\infty}\left(\frac1p\rho_p(u_{n})^p-C_1\rho_p(u_{n})^{\pst}-\|f\|_{\WW_p^*}\rho_p(u_n)\right) \\
 & \geq 2\delta_0-\|f\|_{\WW_p^*} r_0.
\end{align*}
Thus there exists $d_2>0$ such that if $0<\|f\|_{\WW_p^*}<d_2$, it holds $2\delta_0-\|f\|_{\WW_p^*} r_0>\delta_0$. Thus we get  $$0>c_{f,1}>\delta_0>0,$$
which yields a contradiction and hence $u_{n}\in B_{r_0-\eps_0}$ for every $n$. For $0<\eps_1<\eps_0$, we have $B_{r_0-\eps_0}\subset B_{r_0-\eps_1}$. Applying Ekeland's variational principle on the complete metric space $(\ov{B}_{r_0-\eps_1},\rho_p)$, we get
$$c_{f,1}\leq J_f(u_{n})\leq c_{f,1}+\frac{1}{n},$$  $$ J_f(u_{n})\leq J_f(v)+\frac1n\rho_p(u_{n}-v), \; \forall \, v\in \ov{B}_{r_0-\eps_1}, \, v\neq u_{n}.$$
Further, $$J_f(v)=J_f(u_{n})+\prescript{}{\WW_p^*}{\langle}J_f'(u_{n}),(v-u_{n})\rangle_{\WW_p}+o(\rho_p(v-u_{n})).$$
Let $w\in \WW_p$ and $t>0$ be such that $\rho_p(w)=1$ and $v=u_{n}+tw\in \ov{B}_{r_0-\eps_1}$. Therefore, using the above relation we get $$-\frac tn\leq J_f(u_{n}+tw)-J_f(u_{n})=t\prescript{}{\WW_p^*}{\langle}J_f'(u_{n}),w\rangle_{\WW_p}+o_t(1).$$
Dividing by $t$ and letting $t\to0$, we obtain $$-\frac1n\leq \prescript{}{\WW_p^*}{\langle}J_f'(u_{n}),w\rangle_{\WW_p}.$$Replacing $-w$ by $w$, we conclude $$\|J_f'(u_{n})\|_{\WW_p^*}\leq\frac1n.$$
Thus, $\{u_{n}\}$ is a (PS) sequence of $J_f$ at the level $c_{f,1}<0$ in $B_{r_0}$. Suppose $u_n\rightharpoonup v_{f,0}$ in $\WW_p$. Clearly, $v_{f,0}\in B_{r_0}$. By Proposition \ref{grad_p.w.}, we see that for every $\phi\in\WW_p$,
$$\int_{\RR}|\nabla u_n|^{p-2}\nabla u_n\cdot\nabla\phi\dx\to \int_{\RR}|\nabla v_{f,0}|^{p-2}\nabla v_{f,0}\cdot\nabla\phi\dx.$$
For the nonlocal part, we first observe that the sequence $\left\{\frac{A_p(u_n)(x,y)}{|x-y|^{\frac{N+ps}{p'}}}\right\}_n$ is bounded in $L^{p'}(\R^{2N})$. Since $u_n\to v_{f,0}$ a.e. in $\R^N$ (up to a subsequence), it follows that $$\frac{A_p(u_n)(x,y)}{|x-y|^{\frac{N+ps}{p'}}}\to \frac{A_p(v_{f,0})(x,y)}{|x-y|^{\frac{N+ps}{p'}}} \text{ a.e. in }\R^{2N}.$$Therefore,
$$\frac{A_p(u_n)(x,y)}{|x-y|^{\frac{N+ps}{p'}}}\rightharpoonup \frac{A_p(v_{f,0})(x,y)}{|x-y|^{\frac{N+ps}{p'}}}\text{ in }L^{p'}(\R^{2N}).$$ Hence for every $\phi\in\WW_p$,
$$\AA_p(u_n,\phi)\to \AA_p(v_{f,0},\phi).$$ Combining the above convergences and using
 $J_f'(u_n)\to0$, we obtain
\begin{align*}
  \int_{\RR}|\nabla v_{f,0}|^{p-2}\nabla v_{f,0}\cdot\nabla\phi\dx+\AA_p(v_{f,0},\phi) =\int_{\RR}|v_{f,0}|^{p^*-2}v_{f,0}\phi\dx+{}_{\WW_p^*}\langle f,\phi\rangle_{\WW_p}, \; \forall \, \phi \in \mathcal{W}_p.
\end{align*}
Hence $v_{f,0}$ is a solution of \eqref{main_prob_p} in $\WW_p$. Now \eqref{cc-1} and \eqref{cc-6} yield
\begin{align*}
    c_{f,1}&=\limsup_{n\to\infty}\left[ J_f(u_n)-\frac1{p^*}\prescript{}{\WW_p^*}{\langle} J_f'(u_n),u_n\rangle_{\WW_p}\right] \\
    &=\limsup_{n\to\infty}\left[\left(\frac1p-\frac1{p^*}\right)\left(\|\nabla u_n\|_p^p+[u_n]_{s,p}^p\right)-\left(1-\frac1{p^*}\right)\prescript{}{\WW_p^*}{\langle} f,u_n\rangle_{\WW_p}\right]\\
    &\geq \left(\frac1p-\frac1{p^*}\right)\left(\|\nabla v_{f,0}\|_p^p+[v_{f,0}]_{s,p}^p+\mu_{\infty}\right)-\left(1-\frac1{p^*}\right){}_{\WW_p^*}\langle f,v_{f,0}\rangle_{\WW_p}+\frac1N \sum_{j\in J}\mu_j\\
    &\geq J_f(v_{f,0})+ \frac1N \mu_{\infty} + \frac1N S_p^{\frac{N}{p}} \text{card} (J)\geq c_{f,1} + \frac1N \mu_{\infty} +\frac1N S_p^{\frac{N}{p}} \text{card} (J).
\end{align*}
In the penultimate inequality, we use \eqref{mixed3} and the fact that $v_{f,0}$ solves \eqref{main_prob_p}. In the last inequality, we notice that since $v_{f,0}\in B_{r_0}$, $J_f(v_{f,0})\geq c_{f,1}$. Hence we obtain $\text{card} (J)=0$, which gives $J=\emptyset$, and we also obtain $\mu_{\infty} = 0$.
Therefore, in view of the decomposition \eqref{cc-6} and using the uniform convexity of $\WW_p$, we finally obtain $u_n\to v_{f,0}$ in $\WW_p$. Hence $J_f(v_{f,0})=c_{f,1}$ and as a consequence $\abs{v_{f,0}}$ also minimizes $c_{f,1}$. Therefore, without loss of generality, $v_{f,0}$ is nonnegative. Now, using the strong maximum principle in Proposition \ref{SMP}, $v_{f,0}>0$ a.e. in $\RR$. \qed
\section{Nonexistence result for homogeneous problem}\label{non_exst_section}
This section proves the non-existence result for \eqref{main_prob_p} with $f\equiv0$. We first make the following remarks.

\begin{remark}\label{M1}
    In this remark, we observe that every weak solution to \eqref{main_prob_p} with $f \equiv 0$ lies in $\mathcal{C}^{0,\al}_{\text{loc}}(\RR)$ where $\al\in(0,1)$. If $u\in L^{p_s^*}(\RR)$, then using H\"{o}lder's inequality, $$\int_{\RR}\frac{|u(x)|^{p}}{(1+|x|)^{N+sp}}\dx\leq \|u\|_{p_s^*}^{p}\left(\int_{\RR}\frac{\dx}{(1+|x|)^{(N+sp)q}}\right)^{\frac1q}<\infty,$$ where $q$ is the conjugate exponent of $\frac{p_s^*}{p}$. Thus, $$u\in L_{sp}^{p}(\RR):=\left\{u\in L_{\text{loc}}^{p}(\RR):\int_{\RR}\frac{|u(x)|^{p}}{(1+|x|)^{N+sp}}\dx<\infty\right\}.$$ As a consequence, we get $\WW_p \subset L_{sp}^{p}(\RR)$. Further, from Proposition \ref{regularity}, we know that $g= \abs{u}^{p^*-2}u \in L^{\infty}(\R^N)$. Therefore, if $u \in \WW_p$ weakly solves \eqref{main_prob_p} with $f\equiv0$, then $u \in L^{\infty}(\RR) \cap  \WW_p$ weakly solves the following problem
    \begin{equation}\label{reg_eq}
        \del_pu+\frap u=g \text{ in } \Omega, \; g \in L^{\infty}(\Omega),
    \end{equation}
    for every bounded open set $\Omega \subset \RR$. Let $K\subset\RR$ be an arbitrary compact set and choose an open set $\Omega$ such that $K\Subset \Om$. Keeping \cite[Remark 3]{DeMi} in mind, we apply \cite[Theorem 6]{DeMi} to weak solutions of \eqref{reg_eq} to get $u\in \mathcal{C}^{0,\al}_{\text{loc}}(\Om)$ for every $\al\in(0,1)$. Thus, $u\in \mathcal{C}^{0,\al}(K)$ and as a consequence, $u\in \mathcal{C}^{0,\al}_{\text{loc}}(\RR)$ for every $\al\in(0,1)$.
\end{remark}

For $N>p$, consider the problem \begin{align}\label{p1}
 \delp u+\frap u=h(u) \text{ in } \RR,
\end{align}
where $h\in \mathcal{C}(\R)$ given by $h(t)=|t|^{p^*-2}t$. Applying \cite[Corollary 2.7]{AG}, we see that if $u\in \WW_p \cap L^{p}(\RR)\cap L^{\infty}(\RR)\cap \mathcal{C}^{0,l}(\RR)$ weakly solves \eqref{p1}, for some $l \in (s,1)$, then $u$ satisfies the following identity
\begin{equation}\label{p2}
       \frac{N-p}p\|\nabla u\|_p^p+\frac{N-sp}{p}[u]_{s,p}^p=N\int_{\RR} H(u)\dx,
\end{equation}
where $H$ is the primitive of $h$. To use this result for our purpose,  in the following remark we relax the restrictive assumptions on $u$.

\begin{remark}\label{G1} \noi(a) In the proof of \cite[Corollary 2.7]{AG}, we find that the condition $u\in L^{p}(\RR)$ is required only in \cite[pg. 10]{AG} to show that
\begin{equation}\label{AG_cond_1}
    \lim_{t\to0}\|h(u(\cdot+te_j))-h(u)\|_{L^{p'}(\RR)} = 0,
\end{equation}
where $p'$ is the conjugate exponent of $p$ and $e_j$ is the unit vector in the $j$th coordinate. However, we claim that $u\in L^{p}(\RR)$ is not necessary to get \eqref{AG_cond_1}. Indeed, if $u \in \mathcal{W}_p$ weakly solves \eqref{p1}, then applying Proposition \ref{regularity}, we get $u\in L^r(\RR)$ for every $r \in [p^*, \infty]$. Hence using $(p^*-1)p'>p^*$, we see that $h\circ u\in L^{p'}(\RR)$. Now, from the continuity of the translation operator, \eqref{AG_cond_1} holds.

\noi(b) The authors have also used the condition $u\in \mathcal{C}^{0,l}(\RR)$ only in \cite[pg. 11]{AG} to show that
$$I=\lim_{\mu\to0^+}\mu^{-N-1-sp}\iint_{\substack{x,y\in\RR\\|x-y|=\mu}}|u(x)-u(y)|^p(y-x)\cdot\left(\Psi_\la(x)-\Psi_\la(y)\right)\d\si(x,y)=0.$$
Here $\Psi_\la\in \mathcal{C}^{0,1}(\RR,\RR)$ and $\Psi_{\la}$ has compact support. Thus there exists $R>0$ large enough such that $\Psi_\la(x)=\Psi_\la(y)$ for every $x,y\in S_\mu^c$ and for every $\mu\in(0,1)$, where $$S_\mu=\left\{(x,y)\in B_R(0) \times B_R(0):|x-y|=\mu \right\}.$$ Hence, the above integral can be reduced to
$$I=\lim_{\mu\to0^+}\mu^{-N-1-sp}\iint_{S_\mu}|u(x)-u(y)|^p(y-x)\cdot(\Psi_\la(x)-\Psi_\la(y))\d\si(x,y).$$
Now, using $u\in \mathcal{C}^{0,l}(\RR)$ with $l>s$, they get $I=0$. However, since $S_\mu \subset B_R(0) \times B_R(0)$, it is enough to consider $u\in \mathcal{C}_{\text{loc}}^{0,l}(\RR)$ to get $I=0$.

\noi(c)  From (a) and (b), we observe that for a weak solution to \eqref{p1}, Pohozaev's identity \eqref{p2} holds, even when $u\in \WW_p\cap \mathcal{C}_{\text{loc}}^{0,l}(\RR)$.
\end{remark}

In view of Remark \ref{M1} and Remark \ref{G1}, we get the following result.

\begin{proposition}\label{Pohozaev_thm}
   Let $N>p$ and $u\in \WW_p$ be a weak solution to \eqref{main_prob_p} with $f \equiv 0$. Then $u$ satisfies the following identity:
   \begin{equation}
       \frac{N-p}p\|\nabla u\|_p^p+\frac{N-sp}{p}[u]_{s,p}^p=\frac{N}{p^*}\|u\|_{p^*}^{p^*}.
   \end{equation}
\end{proposition}

\noi \textbf{Proof of Theorem \ref{non-existence theorem}:}
    Suppose $u\in \WW_p$ weakly solves \eqref{main_prob_p} with $f\equiv0$. Using $u$ as a test function in \eqref{test-function-p},
    \begin{equation}\label{(a)}
        \|\nabla u\|_p^p+[u]_{s,p}^p=\|u\|_{p^*}^{p^*}.
    \end{equation}
    Applying Proposition \ref{Pohozaev_thm},
    \begin{equation*}
        \frac{N-p}{p} \|\nabla u\|_p^p+\frac{N-sp}{p}[u]_{s,p}^p=\frac{N}{p^*}\|u\|_{p^*}^{p^*}.
    \end{equation*}
    We can rewrite the above equation as
    \begin{equation*}
        \frac{N-p}{p}\left( \|\nabla u\|_p^p+[u]_{s,p}^p\right)+(1-s)[u]_{s,p}^p=\frac{N}{p^*}\|u\|_{p^*}^{p^*}.
    \end{equation*}
    Using \eqref{(a)}, this implies that
    \begin{equation*}
        \frac{N}{p^*}\|u\|_{p^*}^{p^*}+(1-s)[u]_{s,p}^p=\frac{N}{p^*}\|u\|_{p^*}^{p^*}.
    \end{equation*}
    That means $(1-s)[u]_{s,p}^p=0$, which immediately gives  $u=0$ a.e. in $\RR$.
    \qed

\section*{Acknowledgments}	We are grateful to the anonymous reviewers for their very careful corrections and insightful comments. The research of Nirjan Biswas is supported by the National Board for Higher Mathematics Postdoctoral Fellowship (0204/16(9)/2024/RD-II/6761). The research of Paramananda Das is supported by the National Board for Higher Mathematics (0203/5(38)/2024-RD-II/11224).

\bibliographystyle{abbrvnat}

\end{document}